\documentclass[aop,MSNbibl,citesort,seceqn,dvips]{arximspdf}
\usepackage{accents}

\doi{10.1214/11-AOP663}
\volume{40}
\issue{4}
\pubyear{2012}
\firstpage{1795}
\lastpage{1828}

\makeatletter

\newcommand{\ccccdot}{\,\cdot\,}

\newcommand{\cop}[1]{\displaystyle\mathop{#1}^{{\stackrel{\circ}{\mbox{\downparenfill}}}}}
\newcommand{\downparenfill}{$\braceld\rule[0.1pt]{1cm}{0.8pt}\bracerd$}

\newcommand{\R}{\mathbb R}
\newcommand{\N}{\mathbb N}

\newtheorem{Thm}{Theorem}[section]
\newtheorem{Lem}{Lemma}[section]
\newtheorem{Pro}{Proposition}[section]
\newtheorem{Cor}{Corollary}[section]
\newproclaim{Dfn}{Definition}[section]

\newproclaim{Remarks}{Remarks}
\newproclaim{Remark}{Remark}
\newproclaim{AdditionalNotation}{Additional notation}
\newproclaim{Case}{Case}

\newcommand{\stackover}[1]{\accentset{\circ}{#1}}

\newcommand{\ovli}{\overline}

\makeatother

\begin{document}
\begin{frontmatter}

\title{The local quantization behavior of absolutely continuous probabilities}
\runtitle{Local quantization behavior}

\begin{aug}
\author[A]{\fnms{Siegfried} \snm{Graf}\ead[label=e1]{graf@fim.uni-passau.de}},
\author[B]{\fnms{Harald} \snm{Luschgy}\ead[label=e2]{luschgy@uni-trier.de}} and
\author[C]{\fnms{Gilles} \snm{Pag\`{e}s}\corref{}\ead[label=e3]{gilles.pages@upmc.fr}}
\runauthor{S. Graf, H. Luschgy and G. Pag\`{e}s}
\affiliation{Universit\"{a}t Passau, Universit\"{a}t Trier and Universit\'{e} Pierre et Marie Curie}
\address[A]{S. Graf\\
Fakult\"{a}t f\"{u}r Informatik\\
\quad und Mathematik\\
Universit\"{a}t Passau\\
D-94030 Passau\\
Germany\\
\printead{e1}} 
\address[B]{H. Luschgy\\
FB IV-Mathematik\\
Universit\"{a}t Trier\\
D-54286 Trier\\
Germany\\
\printead{e2}}
\address[C]{G. Pag\`{e}s\\
Laboratoire de Probabilit\'{e}s\\
\quad et Mod\`{e}les Al\'{e}atoires\\
UMR~7599, Universit\'{e}\\
\quad Pierre et Marie Curie\\
case 188\\
4, pl. Jussieu\\
F-75252 Paris Cedex 5\\
France\\
\printead{e3}}
\end{aug}

\received{\smonth{10} \syear{2010}}
\revised{\smonth{3} \syear{2011}}

%
\begin{abstract}
For a~large class of absolutely continuous probabilities $P$ it is
shown that, for $r > 0$, for $n$-optimal $L^r(P)$-code\-books
$\alpha_n$, and any Voronoi partition $V_{n,a}$ with respect\vspace*{1pt} to
$\alpha_n$ the local probabilities $P(V_{n,a})$ satisfy $P(V_{a,n})
\approx n^{-1}$ while the local $L^r$-quantization errors
satisfy $\int_{V_{n,a}} \|x-a\|^r \,dP(x) \approx n^{- (1+
r/d)}$ as long as the partition sets~$V_{n,a}$ intersect a~fixed
compact set $K$ in the interior of the support of~$P$.
\end{abstract}

%
\begin{keyword}[class=AMS]
\kwd{60E99}
\kwd{62H30}
\kwd{34A29}.
\end{keyword}
\begin{keyword}
\kwd{Vector quantization}
\kwd{probability of Voronoi cells}
\kwd{inertia of Voronoi cells}.
\end{keyword}

\end{frontmatter}

\section{Introduction}
\label{intro}

The theory of quantization of probability distributions has its origin
in electrical engineering and image processing where it plays a~decisive role in digitizing analog signals and compressing digital
images (see Gray and Neuhoff~\cite{Gray}). More recently, it has also
found many applications in numerical integration (see, e.g.,
\cite{Cher,CherOn,Pages0,Pages}) and mathematical
finance (see, e.g.,~\cite{PAPR} for a~survey).

Optimal (vector) quantization deals with the best approximation of an
$\R^d$-valued random vector $X$ with probability distribution $P$ by
$\R^d$-valued random vectors which attain only finitely many values.
If $r > 0$ and $ \int\|x\|^r \,dP < \infty$ and $n \in\N$, then
the $n$th-level $L^r(P)$-quantization error is defined to be
%
%
\begin{eqnarray}\qquad
\label{eq1}
e_{n,r} &=& e_{n,r}(P) \nonumber\\
&=& \inf\biggl\{ \biggl( \int\|x-q(x)\|^r
\,dP(x) \biggr)^{1/r} \Big\vert q \dvtx \R^d \to\R^d\mbox{ Borel
measurable}\\
&&\hspace*{184pt} \mbox{with }\operatorname{card}(q(\R^d)) \leq n
\biggr\},\nonumber
\end{eqnarray}
where \mbox{$\| \cdot\|$} is a~norm on $\R^d$ and $\operatorname{card}(A)$ stands for
the cardinalility of $A$.

It is known that the above infimum remains unchanged if the Borel
functions $q \dvtx \R^d \to\R^d$ are chosen to be projections onto their
range $\alpha:= q(\R^d) \subset\R^d$ with $\operatorname{card}(\alpha) \leq n$
which obey a~nearest neighbor rule, that is,
\[
q(x) = \sum_{a \in\alpha} a
1_{V_{n,a}}(x),
\]
where $(V_{n,a})_{a \in\alpha}$ is a~Voronoi partition of $\R^d$
with respect to $\alpha$, that is, a~Borel partition such that each of
the partition sets $V_{n,a}$ is contained in the Voronoi cell $W(a
\vert\alpha_n) := \{x \in\R^d \vert\|x-a\| =
\min_{b \in\alpha} \|x-b\|\}$.

If $d(x, \alpha) := \min_{a \in\alpha} \|x-a\|$
denotes the distance of $x$ to the set $\alpha$, then
\[
e_{n,r} = \inf\biggl\{ \biggl( \int d(x, \alpha)^r \,dP(x)
\biggr)^{1/r} \Big\vert\alpha\subset\R^d \mbox{ and } \operatorname{card}(\alpha) \leq n \biggr\}.
\]
The above infimum is in fact a~minimum which is attained at
an optimal ``codebook'' $\alpha_n$ (see~\cite{Graf-Lu}, Theorem
4.12). If $P$ is absolutely continuous with density $h \geq0$ and $
\int\|x\|^{r+\delta} \,dP(x) < \infty$ for some $\delta> 0$, then
%
%
\begin{equation}\label{eq2}
\lim_{n \to\infty} n^{1/d} e_{n,r}(P) = Q_r(P)^{1/r}
\end{equation}
for a~positive real constant $Q_r(P)$ (see Zador~\cite{Zador,Zador-As},
Bucklew and Wise~\cite{Buck} and Graf and Luschgy~\cite{Graf-Lu}, Theorem
6.2). Thus, the sharp asympotics of the sequence $(e_{n,r}^r
)_{n \in\N}$ is completely elucidated up to the numerical
value of the constant $Q_r(P)$.

A famous conjecture of Gersho~\cite{Gersho} states that the bounded
Voronoi-cells of $L^r$-optimal codebooks $\alpha_n$ have
asymptotically the same $L^r$-inertia and a~normalized shape close to
that of a~fixed polyhedron $H$ as $n$ tends to infinity.

In particular, this conjecture suggests that the local
$L^r$-quantization errors (\mbox{$=$}$L^r$-local inertia) satisfy
%
%
\begin{equation}\label{eq3}
\int_{V_{n,a}} \|x-a\|^r \,dP(x) \sim\frac{1}{n}
e_{n,r}^r,\qquad a \in\alpha_n,
\end{equation}
where $a_n \sim b_n$ abbreviates $a_n=\varepsilon_n b_n$ with $\lim
_{n \to\infty} \varepsilon_n = 1$.

So far, this last statement has only been proved for certain parametric
classes of one-dimensional distributions $P$ (see Fort and Pag\`{e}s
\cite{Fort}).

In the present paper, we will investigate the asymptotic behavior for
$n \to\infty$ of $P(W(a \vert\alpha_n))$ and $ \int_{W(a
\vert\alpha_n)} \|x-a\|^r \,dP(x)$ for a~large class of distributions
on $\R^d$ including the nonsingular normal distributions. To derive a~conjecture for the asymptotic size of $P(W(a \vert\alpha_n))$,
one can use the following heuristics. The empirical measure theorem
(see~\cite{Graf-Lu}, Theorem 7.5) states that the empirical
probabilities $\frac{1}{n} \sum_{a \in\alpha_n} \delta
_a$ weakly converge as $n\to\infty$ to the ``point density measure''
%
\[
P_r= \frac{1}{ \int h^{{d}/({r+d})} \,d \lambda^d} h^{{d}/({r+d})}
\lambda^d,
\]
where $\lambda^d$ denotes the $d$-dimensional Lebesgue measure. Thus
we obtain, at least for bounded continuous densities $h$ and an
arbitrary bounded continuous function $f \dvtx\R^d \to\R$, that
%
%
\begin{eqnarray}
\label{eq4}
&&\lim_{n \to\infty} \sum_{a \in
\alpha_n} \frac{1}{n} \biggl( \int h^{{d}/({r+d})} \,d \lambda
^d \biggr) h^{{r}/({r+d})}(a) \int f\,d
\delta_a \nonumber\\
&&\qquad = \int h^{{d}/({r+d})} \,d \lambda^d \lim_{n \to\infty} \biggl(\frac{1}{n}
\sum_{a \in\alpha_n} h^{{r}/({r+d})}(a) f(a) \biggr) \nonumber\\[-8pt]\\[-8pt]
&&\qquad = \int h^{{d}/({r+d})} \,d \lambda^d \int h^{{r}/({r+d})}(x) f(x) \,dP_r(x)\nonumber\\
&&\qquad = \int f(x) \,dP(x),\nonumber
\end{eqnarray}
so that
\[
\sum_{a \in\alpha_n} \biggl( \frac{1}{n} \int h^{{d}/({r+d})} \,d \lambda^d \biggr) h^{{r}/({r+d})} (a)\delta_a
\stackrel{(\R^d)}{\Longrightarrow} P,
\]
where\vspace*{-1pt} $\stackrel{(\R^d)}{\Longrightarrow} $ denotes the weak
convergence of finite measures on $\R^d$. Since it is well known\vspace*{1pt} that
$\sum_{a \in\alpha_n} P(V_{n,a}) \delta_a\stackrel{(\R
^d)}{\Longrightarrow} P$ as well (see~\cite{Pages0,Pages} but also
\cite{Cher,CherOn} or~\cite{Graf-Lu}, Equation (7.6)), it is
reasonable to conjecture that
%
%
\begin{equation}\label{eq5}
P(V_{n,a}) \sim\frac{1}{n} \biggl( \int h^{{d}/({r+d})} \,d \lambda
^d \biggr) h^{{r}/({r+d})}(a).
\end{equation}
We were not able to prove this asymptotical behavior of $P(V_{n,a})$ in
its sharp and general form. But we will show that, for a~large class of
absolutely continuous distributions $P$, there are real constants $c_1,
c_2, c_3, c_4 > 0$ only depending on $P$ such that
$\forall K \subseteq\R^d$, compact, $\exists n_K
\in\N, \forall n \geq n_K, \forall a \in\alpha_n$
%
%
\begin{eqnarray}\label{eq6}
&&K \cap W(a \vert\alpha_n) \neq\varnothing\nonumber\\
&&\quad\Longrightarrow\quad
\frac{c_1}{n} \bigl(\operatorname{essinf} h_{\vert
W_0(a \vert\alpha_n)}\bigr)^{{r}/({r+d})} \\
&&\qquad\qquad\qquad\hspace*{0pt} \le
P(V_{n,a})
\leq\frac{c_2}{n} \bigl(\operatorname{esssup} h_{\vert
W(a \vert\alpha_n)} \bigr)^{{r}/({r+d})},\nonumber
\end{eqnarray}
where
%
%
\begin{equation}\label{W0}
W_0(a \vert\alpha_n) =\bigl\{x \in\R^d \vert\|x-a\| <
d(x, \alpha_n \setminus\{a\})\bigr\}
\end{equation}
and
%
%
\begin{equation}\label{eq7}
\frac{c_3}{n} e^r_{n,r} \leq\int_{V_{n,a}} \|x-a\|^r \,dP(x) \leq
\frac{c_4}{n} e_{n,r}^r.
\end{equation}
The proofs mainly rely on the following two ingredients:

$\bullet$ A ``differentiated Zador's theorem''
%
%
\begin{equation}\label{eq8}
e_{n,r}^r - e^r_{n+1,r} \approx n^{-(1+{r/d} )}
\end{equation}
[where $a_n \approx b_n$ means that the sequence $(\frac
{a_n}{b_n})$ is bounded and bounded away from~$0$] and

$\bullet$ Two \textit{micro--macro inequalities} which relate the
\textit{pointwise distance} of a~quantizer to the global \textit{mean
quantization error} induced on a~distribution~$P$ by this quantizer:

For $b \in(0, \frac{1}{2} )$ fixed, there is
a constant $c_5 > 0$ with
%
%
\begin{equation}\label{eq9}
\quad\qquad\forall n \in\N, \forall x
\in\R^d\qquad c_5 (e^r_{n,r} - e^r_{n+1,r} ) \geq d(x,
\alpha_n)^r P(B(x, bd(x, \alpha_n)))
\end{equation}
and
%
%
\begin{equation}\label{eq10}\qquad
\forall n \geq2\qquad e^r_{n-1,r} - e_{n,r}^r \leq\int
_{V_{n,a}} \bigl(d(x, \alpha_n \setminus\{a\})^r - \|x-a\|^r \bigr) \,dP(x).
\end{equation}
We have stated and established these inequalities in earlier papers:
see especially~\cite{Graf-Pa}; for a~preliminary version of (\ref
{eq10}), see~\cite{Graf-LRat} and for a~one-sided first version
of (\ref{eq8}), see Lemma 3.2 in~\cite{PASA}. They were somewhat
hidden as technical tools inside proofs but their full impact will
become clear here.

The remaining part of the \hyperref[intro]{Introduction} contains a~sketch of the
contents of the paper.
In Section~\ref{2}, we indicate the proofs of the above micro--macro
inequalities and the (weak) asymptotics of quantization error differences.
In Section~\ref{3}, we show that absolutely continuous probabilities $P$ on
$\R^d$, which have a~peakless, connected and compact support as well
as a~density which is bounded and bounded away from $0$ on the support,
have asymptotically uniform local quantization errors (Theorem~\ref{39}).
In Section~\ref{4}, we show that absolutely continuous probabilities whose
densities are the composition of a~decreasing function on $\R_+$ and a~norm or a~quasi-concave function outside a~compact set satisfy a~sharpened first micro--macro inequality of the following type:

There exist a~constant $c > 0$ such that, for every $K \subset\R^d$ compact,
\[
\exists n_K \in\N, \forall n \geq n_K, \forall x \in K
\qquad c n^{-1/d} h(x)^{-{1}/({r+d})} \geq d(x, \alpha_n).
\]

Assuming this inequality, we derive asymptotic estimates for the
probabilities of the quantization cells and local quantization errors
(Theorem~\ref{49}).
Section~\ref{sec5} deals with the local quantization behavior of certain Borel
probabilities $P$ in the interior of their support. The results are
stated for arbitrary absolutely continuous probabilities with density
$h$ satisfying the moment condition $ \int\|x\|^{r+ \delta}
h(x) \,d \lambda(x) <+ \infty$ for some $\delta> 0$. They are
particularly useful if the density $h$ is bounded and bounded away from
$0$ on each compact subset of the interior of the support of~$P$. Under
these very general assumptions,
the results are quite similar to those given in Section~\ref{4} but the given
constants are a~little bit less effective (Theorem~\ref{56}).
\begin{AdditionalNotation*} For $x \in\R^d$ and
$\rho> 0$ $B(x,\rho) =B_{{\|\ccccdot\|}}(x,\rho)= \{y \in\R^d
\vert\|y-x\| < \rho\}$ denotes the open ball with center $x$ and
radius $\rho$. \mbox{$\| \cdot\|_2$} will denote the canonical Euclidean norm on
$\R^d$.\vspace*{2pt}

$\stackover{A}$ denotes the interior of a~set $A \subset\R^d$.
\end{AdditionalNotation*}

\section{Important inequalities in quantization}
\label{2}

In the following, \mbox{$\| \cdot\|$} denotes an arbitrary norm on $\R^d$ and
$P$ is always an absolutely continuous Borel probability on $\R^d$
which has density $h$ with respect to the $d$-dimensional Lebesgue
measure $\lambda^d$. Let $r \in(0,+ \infty)$ be fixed. We always
assume that there is a~$\delta> 0$ with $ \int\|x\|^{r+ \delta}
\,dP(x) < + \infty$. For every $n \in\N$, let $e_{n,r}$ denote the
$n$th-level $L^r(P)$-quantization error. Then we have
%
%
\begin{equation}\label{eq12}\quad
e_{n,r}^r = e^r_{n,r}(P) = \inf\biggl\{ \int d(x,\alpha)^r \,dP(x)
\Big\vert\alpha\subset\R^d, \operatorname{card}(\alpha) \leq n
\biggr\}.
\end{equation}
For each $n \in\N$, we choose an arbitrary $n$-optimal set $\alpha
_n \subset\R^d$, that is, a~set $\alpha_n \subset\R^d$ with $\operatorname{card}(\alpha_n) \leq n$ and
%
%
\begin{equation}\label{eq13}
e^r_{n,r} = \int d(x, \alpha_n)^r \,dP(x).
\end{equation}
It is well known that, under the above conditions, such a~set exists
and satisfies
%
%
\begin{equation}\label{eq14}
\operatorname{card}(\alpha_n) = n.
\end{equation}
In this section, we will state the fundamental inequalities which
relate the behavior of the distance function $d(\cdot, \alpha_n)$ to
the difference $e^r_{n,r} - e^r_{n+1,r}$ of successive $r$th powers of
the quantization errors. Using these inequalities, we will be able to
determine the (weak) asymptotics of $e^r_{n,r} - e^r_{n+1,r}$.

\subsection{Micro--macro inequalities}
\begin{Pro}[(First micro--macro inequality)] \label{21} For every $b
\in(0, \frac{1}{2} )$, for all $n \in\N$ and all $x
\in\R^d$,
%
%
\begin{equation}\label{eq15}
e_{n,r}^r - e_{n+1,r}^r \geq(2^{-r} - b^r ) d(x, \alpha
_n)^r P(B(x, bd(x, \alpha_n))).
\end{equation}
\end{Pro}
\begin{pf}
The proof can be found as part of the proof of
Theorem 2 in~\cite{Graf-Pa}.
\end{pf}
\begin{Remarks*} (a) Inequality (\ref{eq15}) holds for
arbitrary Borel probabilities $P$ on $\R^d$ for which $ \int\|x\|
^r \,dP(x) < \infty$. $P$ need not be absolutely continuous.

(b) By the differentiation theorem for absolutely continuous
measures $P= h \lambda^d$ and the fact (see~\cite{Delattre}) that
$\lim_{n \to\infty} d(x, \alpha_n) = 0$ for every $x \in
\operatorname{supp}(P)$, we know that, for $\lambda^d$-a.e. $x \in\R^d$,
%
%
\begin{equation}\label{eq16} \lim_{n \to\infty} \frac
{P(B(x,bd(x, \alpha_n)))}{\lambda^d(B(x, bd(x, \alpha_n)))} = h(x).
\end{equation}
\end{Remarks*}

Having this in mind, we can rephrase (\ref{eq15}) as follows:
%
%
\begin{eqnarray}\label{eq17}
&&\forall n \in\N, \forall x \in\R^d\nonumber\\[-8pt]\\[-8pt]
&&\qquad
c_5 (e_{n,r}^r - e^r_{n+1,r}) \geq d(x, \alpha_n)^{r+d} \frac
{P(B(x, bd(x, \alpha_n)))}{\lambda^d(B(x, bd(x, \alpha_n)))} ,\nonumber
\end{eqnarray}
where
\[
c_5=[ (2^{-r}-b^r ) b^d \lambda^d(B(0,1))
]^{-1}
\]
(with the convention $0 \cdot\mbox{undefined}= 0$).

Suppose that there is a~constant $c_9 > 0$ such that
%
%
\begin{equation}\label{eq18}\quad
\exists n_0 \in\N, \forall n \geq n_0, \forall x \in\R
^d\qquad \frac{P(B(x, bd(x, \alpha_n)))}{\lambda^d(B(x, bd(x, \alpha
_n)))} \geq c_9 h(x).
\end{equation}
Then, for $c_{10} = c_5 c_9^{-1}$, we have
%
%
\begin{equation}\label{eq19}
\forall n \geq n_0, \forall x \in\R^d\qquad c_{10}
(e^r_{n,r} - e^r_{n+1,r} ) \geq d(x, \alpha_n)^{r+d} h(x).
\end{equation}

\begin{Pro}[(Second micro--macro inequality)] \label{23} One has
%
%
\begin{eqnarray}\label{eq20}\quad
&&\forall n \geq2, \forall a \in\alpha_n \nonumber\\[-8pt]\\[-8pt]
&&\qquad e^r_{n-1,r} -
e^r_{n,r} \leq\int_{W_0(a \vert\alpha_n)}
\bigl(d(x, \alpha_n \setminus\{a\})^r - \|x-a\|^r \bigr)
\,dP(x),\nonumber
\end{eqnarray}
where $W_0(a \vert\alpha_n)$ is defined by (\ref{W0}).
\end{Pro}
\begin{pf}
The proof is part of the proof of
\cite{Graf-Pa}, Theorem 2.
\end{pf}
\begin{Remark*} Inequality (\ref{eq20}) holds for arbitrary
Borel probabilities $P$ on~$\R^d$ with $\int\|x\|^r \,dP(x) < +\infty$.
\end{Remark*}

\subsection{A differentiated version of Zador's theorem}
To use the preceding propositions for concrete calculations, it is
essential to know the asymptotic behavior of the error differences
$e^r_{n,r} - e^r_{n+1,r}$. We have the following result in that direction.
\begin{Pro}\label{25}
If $P$ is absolutely continuous on $\R^d$, then
\[
e^r_{n,r} - e^r_{n+1,r} \approx n^{- (1 + {r/d}
)}.
\]
\end{Pro}
\begin{pf}
In the proof of Theorem 2 in~\cite{Graf-Pa}, it
is shown that there is a~constant $c_{11} > 0$ such that
\[
\forall n \in\N\qquad e_{n,r}^r - e^r_{n+1,r} \leq
c_{11}n^{-( 1+r/d )}.
\]

To obtain the lower bound for $e^r_{n,r} - e^r_{n+1,r}$, we proceed as
follows.\vadjust{\goodbreak}

It follows from (\ref{eq16}) and Egorov's theorem (see~\cite{Cohn},
Proposition 3.1.3) that there exists a~real constant $c > 0$ and a~Borel set $A \subset\{h>c\}$ of finite and positive Lebesgue measure
such that
\[
\mbox{the convergence of } \frac{P(B(x,bd(x, \alpha_n)))}{\lambda
^d(B(x,bd(x, \alpha_n)))} \mbox{ to } h \qquad\mbox{is uniform in } x
\in A.
\]
Hence, there exists an $n_0 \in\N$ with
%
%
\begin{equation}\label{eq22}
\forall n \geq n_0 , \forall x \in A \qquad \frac{P(B(x,bd(x,
\alpha_n)))}{\lambda^d(B(x,bd(x, \alpha_n)))} > \frac{1}{2} c.
\end{equation}
Combining (\ref{eq17}) and (\ref{eq22}) and integrating both sides of
the resulting inequality with respect to the Lebesgue measure on $A$ yields
\begin{eqnarray*}
c_5 (e^r_{n,r} - e^r_{n+1,r} ) & \geq &\frac{1}{\lambda^d
(A)} \frac{1}{2} c \int_A d(x, \alpha_n)^{r+d} \,d \lambda^d(x)\\
& \geq &\frac{1}{2} c e^{r+d}_{n,r+d} (\lambda^d (\cdot
\vert A)),
\end{eqnarray*}
where $\lambda^d(\cdot\vert A)$ denotes the normalized Lebesgue
measure on $A$.
By Zador's theorem (see (\ref{eq2}) or~\cite{Graf-Lu}, Theorem 6.2),
we have
\[
\liminf_{n \to\infty} n^{1+{r/d}}
e^{r+d}_{n,r+d} (\lambda^d( \cdot|A)) > 0,
\]
so that $\liminf_{n \to\infty} n^{1+{r/d}}
(e^r_{n,r} - e^r_{n+1,r} )> 0$.
\end{pf}
%
\begin{Remark*}
It would be interesting to know the sharp asymptotic behavior of
$e_{n,r}^r - e_{n+1,r}^r$. We conjecture that
\[
\lim_{n \to\infty} n^{1+r/d} (e^r_{n,r} -
e^r_{n+1,r} ) = \frac{d}{r} Q_r(P) = \frac{d}{r} Q_r
([0,1]^d ) \|h\|_{{d}/({d+r})},
\]
where $Q_r([0,1]^d ) \in(0,\infty)$ is as in~\cite{Graf-Lu}, Theorem 6.2.
\end{Remark*}
%
\section{Uniform local quantization rate for absolutely continuous
distributions with peakless connected compact support}
\label{3}

As before, $P$ is an absolutely continuous probability with density
$h$. Let $(\alpha_n)_{n \in\N}$ be a~sequence of optimal codebooks
of order $r \in(0, \infty)$ for $P$. We will investigate the
asymptotic size of
\[
W(a \vert\alpha_n),\qquad P(W(a \vert\alpha_n))
\quad\mbox{and}\quad \int_{W(a \vert\alpha_n)} \|x-a\|^r \,dP(x)
\]
under some compactness and regularity assumptions on $\operatorname{supp}(P)$
and $P$.

The main result of this section is stated below. Its proof, which
heavily relies on the following two paragraphs devoted to upper and
lower bounds, respectively, is postponed to the end of this section.
\begin{Thm} \label{39}Suppose that $P$ is an absolutely continuous
Borel probability on $\R^d$ whose density is essentially bounded,
whose support is connected and compact, and which is ``peakless'' in the
following sense:
\[
\exists c > 0, \exists s_0 > 0, \forall s \in(0,s_0),
\forall x \in\operatorname{supp}(P) \qquad P(B(x,s)) \geq c \lambda^d(B(x,s)).
\]
Let $(\alpha_n)$ be a~sequence of codebooks which are optimal of order
$r \in(0, \infty)$. For $a \in\alpha_n$, let us define the
inradius and the circumradius of the Voronoi cell $W(a \vert
\alpha_n)$ by
\[
\underline{s}_{n,a} = \sup\{s > 0 , B(a,s) \subset W(a \vert
\alpha_n)\}
\]
and
\[
\overline{s}_{n,a} = \inf\{s > 0 , W(a \vert\alpha_n) \cap
\operatorname{supp}(P) \subset B(a,s)\},
\]
respectively. Then
%
%
\begin{eqnarray}
\label{eq44}
\frac{1}{n} &\preccurlyeq&\min_{a \in\alpha
_n} P(W_0(a \vert\alpha_n)) \leq\max_{a \in\alpha_n} P(W(a \vert
\alpha_n)) \preccurlyeq
\frac{1}{n},
\\
\label{eq45} \frac{ e^r_{n,r} }{n} &\preccurlyeq&
\min_{a \in\alpha_n}\int_{W_0(a \vert\alpha_n)}
\|x-a\|^r \,dP(x) \nonumber\\[-8pt]\\[-8pt]
&\leq&\max_{a \in\alpha_n}\int_{W(a
\vert\alpha_n)} \|x-a\|^r \,dP(x)
\preccurlyeq\frac{e^r_{n,r}}{n}\nonumber
\end{eqnarray}
and
%
%
\begin{equation}\label{eq46}
n^{-1/d}\preccurlyeq\min_{a \in\alpha_n} \underline
{s}_{n,a} \leq\max_{a \in\alpha_n} \overline{s}_{n,a}
\preccurlyeq
n^{-1/d}.
\end{equation}
[Here $a_n \preccurlyeq b_n$ means that $( \frac{a_n}{b_n}
)$ is bounded from above.]
\end{Thm}
\begin{Remarks*}
The name ``peakless'' given to the
above assumption illustrates that a~subset of $\R^d$ that satisfies
this condition cannot have infinitely thin peaks (or spines) on its
boundary and that the existence of such peaks or spine is the only way
to make the assumption fail.

Inequality (\ref{eq46}) was proved by Gruber
in~\cite{Gruber}, Theorem 3(ii), under an additional continuity
assumption on $h$, but with a~more general distortion measure.
\end{Remarks*}

\subsection{Upper bounds}
The following proposition is essentially contained in
Graf and Luschgy~\cite{Graf-LRat} (Proposition 3.3 and the following
remark). It has been independently proved by Gruber~\cite{Gruber},
Theorem 3(ii).
\begin{Pro} \label{31}
Suppose that $\operatorname{supp}(P)$ is compact and that there exist constants
$c_{12} > 0$ and $s_0 > 0$ such that
%
%
\begin{equation}\label{eq23}\quad
\forall s \in(0,s_0), \forall x \in\operatorname{supp}(P) \qquad
P(B(x,s))\geq c_{12} \lambda^d(B(x,s)).
\end{equation}
Then there is a~constant $c_{13} < + \infty$ such that
%
%
\begin{equation}\label{eq24}
\forall n \in\N, \forall x \in\operatorname{supp}(P) \qquad
d(x, \alpha_n) \leq c_{13} n^{-1/d}.
\end{equation}
\end{Pro}
\begin{pf}
Let $b \in(0, \frac{1}{2} )$ be fixed.
%
Since $K := \operatorname{supp}(P)$ is compact it follows from~\cite{Delattre},
Proposition 1, that $\lim_{n \to\infty} \max_{x
\in K} d(x, \alpha_n) = 0$. Thus, there is an $n_0 \in\N$ with
\[
\forall n \geq n_0, \forall x \in K \qquad d(x, \alpha_n)
< s_0
\]
and, hence, by (\ref{eq23})
%
%
\begin{equation}\label{eq26}\qquad
\forall n \geq n_0, \forall x \in K\qquad P(B(x, bd(x,
\alpha_n))) \geq c_{12} \lambda^d(B(x,bd(x, \alpha_n))).
\end{equation}
By Proposition~\ref{25}, there exists a~constant $c_{11} > 0$
such that
%
%
\begin{equation}\label{eq27}
\forall n \in\N\qquad e^r_{n,r} - e^r_{n+1,r} \leq c_{11}n^{-
(1+ r/d)}.
\end{equation}
Combining (\ref{eq17}), (\ref{eq26}) and (\ref{eq27}), yields
\[
c_{12}^{-1} c_{11} c_5 n^{- (1+r/d)} \geq d(x, \alpha_n)^{r+d}
\]
for every $x \in K$ and every $n \geq n_0$. Inequality (\ref{eq24})
follows by setting
\[
c_{13} = \max\bigl\{(c^{-1}_{12} c_{11} c_5 )^{
{1}/({r+d})}, \max\bigl\{d(x, \alpha_n ) n^{1/d} , x \in K, n
\in\{1, \ldots, n_0 \}\bigr\} \bigr\}.\quad
\]
\upqed\end{pf}
\begin{Pro}[(Upper-bounds)] \label{32}Suppose that the assumptions of
Pro\-position~\ref{31} are satisfied and that, in addition, $h$ is
essentially bounded. Then there exist constants $c_{14}, c_{15} \in
(0,\infty)$ such that
%
%
\begin{equation}\label{eq28}\quad
\forall n \in\N, \forall a \in\alpha_n \qquad \cases{
\displaystyle P(W(a \vert\alpha_n)) \leq\frac
{c_{14}}{n},\vspace*{2pt}\cr
\displaystyle \int_{W(a \vert\alpha_n)} \|x-a\|^r \,dP(x) \leq
c_{15} n^{- (1+r/d)}.}
\end{equation}
\end{Pro}
\begin{pf}
By Proposition~\ref{31}, we have, for every $n
\in\N$ and every $a \in\alpha_n$,
\[
W(a \vert\alpha_n) \cap\operatorname{supp}(P) = \{x
\in\operatorname{supp}(P) \vert\|x-a\| = d(x, \alpha_n)\} \subseteq B
(a,c_{13} {n}^{-1/d} ),
\]
which implies
\begin{eqnarray*}
P(W(a \vert\alpha_n)) & \leq & P(B (a,c_{13}
{n}^{-1/d} )
) = \int_{B (a,c_{13} {n}^{-1/d} )}
h\,d \lambda^d\\
& \leq & \|h\|_{\R^d} \lambda^d(B(0,1)) c_{13}^d \frac{1}{n},
\end{eqnarray*}
where
$\|h\|_B = \operatorname{esssup} h_{|B}$.
Likewise, we obtain
\begin{eqnarray*}
\int_{W(a \vert\alpha_n)} \|x-a\|^r \,dP(x) & \leq & \int_{B
(a,c_{13} {n}^{-1/d} ) } \|x-a\|^r \,dP(x)\\
& \leq &(c_{13} n^{-1/d} )^r P (B (a,c_{13}
{n}^{-1/d} ) ).
\end{eqnarray*}
Setting $ c_{14} = \|h\|_{\R^d} \lambda
^d(B(0,1)) c_{13}^d$ and
$c_{15} = c_{14} c_{13}^{r}$
yields (\ref{eq28}).
\end{pf}
\begin{Remark*}
Thus, assumption (\ref{eq23}) is satisfied if
$\operatorname{supp}(P)$ is peakless, that is,
%
%
\begin{eqnarray}\label{eq29}
&&\exists c > 0, \exists s_1 > 0, \forall s \in(0,s_1),
\forall x \in\operatorname{supp}(P)
\nonumber\\[-8pt]\\[-8pt]
&&\qquad\lambda^d\bigl(B(x,s) \cap\operatorname{supp}(P)\bigr) \geq c
\lambda^d(B(x,s)),\nonumber
\end{eqnarray}
and $h$ is essentially bounded away from $0$ on $\operatorname{supp}(P)$, that
is,
\[
\exists\underline t > 0 , h(x) \geq\underline t \qquad
\mbox{for } \lambda^d\mbox{-a.e. }x \in\operatorname{supp}(P).
\]
As an example, (\ref{eq29}) holds for finite unions of compact convex
sets with positive $\lambda^d$-measure (see~\cite{Graf-Lu}, Example
12.7 and Lemma 12.4).
\end{Remark*}

\subsection{Lower bounds}
\begin{Lem}\label{34} If $\operatorname{supp}(P)$ is connected then, for every
$n\ge2$ and every $a \in\alpha_n$,
%
%
\begin{equation}\label{eq30}
d(a, \alpha_n \setminus\{a\}) \leq2 \sup\bigl(\{\|y-a\| , y \in
W(a \vert\alpha_n) \cap\operatorname{supp}(P)\}\bigr).
\end{equation}
\end{Lem}
\begin{pf}
Let $n \geq2$ be fixed. First, we will show that
%
%
\begin{equation}\label{eq31}
\forall a \in\alpha_n\qquad W(a \vert\alpha_n) \cap
\bigcup_{b \in\alpha_n \setminus\{a\}}
W(b \vert\alpha_n) \cap\operatorname{supp}(P) \neq
\varnothing.
\end{equation}
Let $a \in\alpha_n$. Since the nonempty closed sets (see \cite
{Graf-Lu}, Theorem 4.1) $W(a \vert\alpha_n) \cap\operatorname{supp}(P)$ and $\bigcup_{b \in\alpha_n \setminus\{a\}} W(b
\vert\alpha_n) \cap\operatorname{supp}(P)$ cover the connected set
$\operatorname{supp}(P)$, claim (\ref{eq31}) follows.

By (\ref{eq31}), there exists $b \in\alpha_n \setminus\{
a\}$ with \mbox{$W(a \vert\alpha_n) \cap W(b \vert\alpha_n)
\cap\operatorname{supp}(P)\neq\varnothing$}.

Let $z$ be a~point in this set. Then $\|z-a\| = d(z,
\alpha_n) = \|z-b\|$ and
\begin{eqnarray*}
d(a, \alpha_n \setminus\{a\}) & \leq &\|a-b\| \leq\|a-z\| + \|z-b\|\\
& \leq &2\|z-a\| \leq2 \sup\{\|y-a\| , y \in W(a \vert\alpha
_n) \cap\operatorname{supp}(P)\}.
\end{eqnarray*}
\upqed\end{pf}
%
\begin{Pro}[(Lower bounds I)] \label{35} Suppose that $\operatorname{supp}(P)$ is
compact and connected, that $P$ satisfies (\ref{eq23}) and is
absolutely continuous with an essentially bounded probability density $h$.

Then there exist constants $c_{16}, c_{17} > 0$ such that
%
%
\begin{equation}\label{eq32a}
\forall n \geq2, \forall a \in\alpha_n \qquad d(a,
\alpha_n \setminus\{a\}) \geq c_{16} n^{-1/d}
\end{equation}
and
%
%
\begin{equation}\label{eq32b}
\forall n \in\N, \forall a \in\alpha_n \qquad P(W_0(a
\vert\alpha_n)) \geq\frac{c_{17}}{n}.
\end{equation}
\end{Pro}
\begin{pf}
Let $n \geq2$ and $a \in\alpha_n$ be
arbitrary. By the second micro--macro inequality (\ref{eq20}), we have
%
%
\begin{eqnarray}\label{eq33}\quad
&&
e^r_{n-1,r} - e^r_{n,r} \nonumber\\
&&\qquad \leq \int_{W_0(a \vert
\alpha_n)} \bigl(d(x, \alpha_n \setminus\{a\})^r - \|
x-a\|^r\bigr) \,dP(x)\\
&&\qquad \leq \int_{W_0(a \vert\alpha_n)}
\bigl( \bigl(\|x-a\| +d(a, \alpha_n \setminus\{a\})\bigr)^r -\|x-a\|^r
\bigr) \,dP(x).\nonumber
\end{eqnarray}
By Proposition~\ref{25}, there exists a~real constant $c >
0$ with
%
%
\begin{equation}\label{eq34} c n^{- (1+r/d)} \leq e^r_{n-1,r} -
e^r_{n,r}.
\end{equation}

\begin{Case}[($r \geq1$)]\label{case1} Combining (\ref{eq33}) and (\ref
{eq34}) and using the mean value theorem for differentiation yields
%
%
\begin{eqnarray}\label{eq35}
c n^{- (1+r/d)} &\leq&
\int_{W_0(a \vert\alpha_n)} r\bigl(\|x-a\| +
d(a, \alpha_n \setminus\{a\})\bigr)^{r-1} \nonumber\\[-8pt]\\[-8pt]
&&\hspace*{37pt}{}\times d(a, \alpha_n \setminus\{a\}
) \,dP(x).\nonumber
\end{eqnarray}
Using Lemma~\ref{34} and (\ref{eq24}), we know that
%
%
\begin{equation}\label{eq36}
\forall x \in W(a \vert\alpha_n) \cap\operatorname{supp}(P)
\qquad \|x- a\| + d(a, \alpha_n \setminus\{a\}) \leq3c_{13} n^{-1/d}.\hspace*{-35pt}
\end{equation}
Combining (\ref{eq35}) and (\ref{eq36}) yields
%
%
\begin{equation}\label{eq37}
r^{-1}c (3c_{13} )^{-(r-1)} n^{-1-1/d} \leq d(a, \alpha_n
\setminus\{a\}) P(W_0(a \vert\alpha_n)).
\end{equation}
Since $P(W_0 (a \vert\alpha_n)) \leq P(W(a \vert
\alpha_n)) \leq c_{14} n^{-1}$ by (\ref{eq28}), we deduce
\[
c^{-1}_{14} r^{-1} c(3c_{13})^{-(r-1)} n^{-1/d} \leq d(a, \alpha
_n \setminus\{a\})
\]
and, hence, (\ref{eq32a}) with $c_{16} = c_{14}^{-1} r^{-1}
c(3c_{13})^{-(r-1)}$.

Since $d(a, \alpha_n \setminus\{a\}) \leq2c_{13} n^{-1/d}$, we
deduce from (\ref{eq37}) that
\[
(2c_{13})^{-1} r^{-1} c(3c_{13})^{-(r-1)} n^{-1} \leq P(W_0(a
\vert\alpha_n))
\]
and, hence, (\ref{eq32b}) with $c_{17} = (2c_{13})^{-1} r^{-1}
c(3c_{13})^{-(r-1)}$.
\end{Case}
\begin{Case}[($r < 1$)]\label{case2}
In this case, we have
\[
\bigl(\|x-a\| +d(a, \alpha_n \setminus\{a\})\bigr)^r \leq\|x-a\|^r +d(a,
\alpha_n \setminus\{a\})^r
\]
for all $x \in W_0(a \vert\alpha_n)$. Combining this
inequality with (\ref{eq33}) and (\ref{eq34}) yields
\[
c n^{- (1+r/d)}
\leq d(a, \alpha_n \setminus\{a\})^r P(W_0(a \vert\alpha_n)).
\]
Since $P(W_0(a \vert\alpha_n)) \leq c_{14}/n$
by (\ref{eq28}), we deduce
\[
(c_{14}^{-1} c )^{1/r} n^{-1/d} \leq d(a, \alpha_n
\setminus\{a\})
\]
and hence, (\ref{eq32a}) with $c_{16} = (c_{14}^{-1} c )^{1/r}$.\vadjust{\goodbreak}

Since $d(a, \alpha_n \setminus\{a\})^r \leq(3c_{13})^r n^{-r/d}$,
we obtain
\[
(3c_{13} )^{-r} c n^{-1} \leq P(W_0(a \vert
\alpha_n) )
\]
and, hence, (\ref{eq32b}) with $c_{17} = (3c_{13} )^{-r}
c$.\qed
\end{Case}
\noqed\end{pf}
\begin{Cor} \label{36} Let the assumptions of Proposition~\ref{35} be
satisfied.

Then there exists a~constant $c_{18} > 0$ such that
%
%
\begin{equation}\label{eq38}
\forall n \in\N, \forall a \in\alpha_n \qquad B
(a, c_{18} n^{-1/d} ) \subset W_0(a \vert\alpha_n).
\end{equation}
\end{Cor}
\begin{pf}
Set $c_{18} = \frac{1}{2} c_{16}$. For
$n=1$ and $a \in\alpha_n$, the assertion is obviously true since
$W_0(a \vert\alpha_1) = \R^d$. Now let $n \geq2$ and let $a
\in\alpha_n$ be arbitrary. We will show that
\[
B (a, c_{18} n^{-1/d} ) \subset W_0(a \vert\alpha_n).
\]
Let $x \in\R^d$ with $\|x-a\| < c_{18} n^{-1/d}$. By
(\ref{eq32a}), we know that
\[
\|x-a\| < \tfrac{1}{2} d(a, \alpha_n \setminus\{a\})
\]
and, hence, for every $b \in\alpha_n \setminus\{a\}$:
\begin{eqnarray*}
\|x-b\| &\geq&\|a-b\|- \|x-a\|  \\
&\geq& d(a, \alpha_n \setminus\{a\}) -
\|x-a\|
> \tfrac{1}{2} d(a, \alpha_n \setminus\{a\})\\
& > &\|x-a\|.
\end{eqnarray*}
This implies $x \in W_0(a \vert\alpha_n)$.
\end{pf}
%
%
%
\begin{Pro}[(Lower bounds II)]\label{38} Let the assumptions of
Proposition~\ref{35} be satisfied. Then there exists a~real constant
$c_{19} > 0$ such that
%
%
\begin{equation}\label{eq39}\quad
\forall n \in\N, \forall a \in\alpha_n \qquad \int
_{W_0(a \vert\alpha_n)} \|x-a\|^r \,dP(x) \geq
c_{19} n^{-(1+{r/d} )}.
\end{equation}
\end{Pro}
\begin{pf}
Let $n \in\N$ and $a \in\alpha_n$ be
arbitrary. By (\ref{eq32b}), we have\break $P(W_0(a \vert\alpha_n)) >
0$. Let $s_a = \inf\{s > 0 \vert P(B(a,s)) \geq\frac
{1}{2} P(W_0(a \vert\alpha_n))\}$. Since $s \mapsto P(B(a,s))$
is continuous with $\lim_{s \downarrow0} P(B(a,s)) = 0$
and\break $\lim_{s \uparrow+ \infty} P(B(a,s)) = 1$, we deduce
%
%
\begin{equation}\label{eq40}
P(B(a,s_a)) = \tfrac{1}{2} P(W_0(a \vert\alpha_n)).
\end{equation}
This implies
%
%
\begin{eqnarray}
\label{eq41}\int_{W_0(a \vert\alpha_n)} \|x-a\|^r \,dP(x)
& \geq &
\int_{W_0(a \vert\alpha_n) \setminus
B(a,s_a)}
\|x-a\|^r \,dP(x) \nonumber\\
& \geq &s_a^r P\bigl(W_0(a \vert\alpha_n) \setminus
B(a,s_a)\bigr)\nonumber\\[-8pt]\\[-8pt]
& \geq &s_a^r\bigl(P(W_0(a \vert\alpha_n)) - P(B(a, s_a))\bigr)\nonumber\\
&=&\frac{1}{2} s^r_a P(W_0(a \vert\alpha_n)).\nonumber
\end{eqnarray}
On the other hand, since $h$ is essentially bounded we have
%
\begin{eqnarray*}
P(W_0(a \vert\alpha_n)) & = & 2P(B(a,s_a))\\
& \leq & 2 \lambda^d(B(a,s_a)) \|h\|_{\R^d}\\
& = & 2 \lambda^d(B(0,1)) s^d_a \|h\|_{\R^d}.
\end{eqnarray*}
Hence,
%
%
\begin{equation}\label{eq42}
s^r_a \geq\biggl( \frac{1}{2 \lambda^d(B(0,1)) \|h\|_{\R^d}}
\biggr)^{r/d} P(W_0(a \vert\alpha_n))^{r/d}.
\end{equation}
Setting $c = \frac{1}{2} ( \frac{1}{2 \lambda
^d(B(0,1)) \|h\|_{\R^d}} )^{r/d}$ and combining (\ref
{eq41}) and (\ref{eq42}) yields
%
%
\begin{equation}\label{eq43}
\int_{W_0(a \vert\alpha_n)}
\|x-a\|^r \,dP(x)
\geq c P(W_0(a \vert\alpha_n))^{1+r/d}.
\end{equation}
Since $P(W_0(a \vert\alpha_n)) \geq c_{17} \frac
{1}{n}$ by (\ref{eq32a}), we deduce
\[
\int_{W_0(a \vert\alpha_n)}
\|x-a\|^r \,dP(x)
\geq c c_{17}^{1+r/d} n^{- (1+r/d)}
\]
and, hence, the conclusion (\ref{eq39}) of the proposition with
$c_{19} = c c_{17}^{1+r/d}$.
\end{pf}
\begin{pf*}{Proof of Theorem~\ref{39}} The result\vspace*{1pt} is a~combination
of the results in Propositions~\ref{31}--\ref
{38}, Corollary~\ref{36} and Zador's theorem which says that $\lim
_{n \to\infty} \frac{e^r_{n,r}}{n^{-r/d}}$ exists in
$(0,+ \infty)$ (see, e.g.,~\cite{Graf-Lu}, Theorem 6.2).
\end{pf*}

\section{The local quantization rate for a~class of absolutely
continuous probabilities with unbounded support}
\label{4}

In this section, we propose extensions of the results of Section \ref
{3} to distributions with an unbounded support which requires to have a~control of the behavior of the distribution at infinity, even if our
results are only \textit{locally} uniform.

First, we introduce in item (c) of the definition below a~class of
probability density functions satisfying the ``Peakless Sublevel Tail
Property'' (PSTP) for which a~sharpened version of the \textit{micro--macro}
inequality (\ref{eq17}) holds [see~(\ref{eq51}) further on]. This
improved inequality is in fact the key to get the main results of this
section (Proposition~\ref{48} and Theorem~\ref{41}).

Although the PSTP may look rather technical and will not be shown to be
necessary for the results in the unbounded framework, it seems clear
from the case of compactly supported distribution, that one needs a~restrictive condition of this nature for the conclusions in the case of
distributions with unbounded support. The ``Peakless Sublevel
Property''
(PSP) [item (a) in the definition below] is in some way the ``core''
of the PSTP and the ``Convex Sublevel Approximation Property'' (CSAP)
[item (b) in the definition below] is simply a~tractable criterion
for the PSP.
\begin{Dfn}\label{41} 
\begin{enumerate}[(a)]
\item[(a)]
A Borel measurable map $f \dvtx\R^d \to\R$ satisfies the
\textit{peakless sublevel property} (PSP) outside $\overline{B}(0,R)$,
$R>0$, if there are real constants $s_0, c_f > 0$ such that
%
%
\begin{eqnarray}
\label{eq47} &&\forall x \in\R^d \setminus\overline
{B}(0,R), \forall s \in(0, s_0)\nonumber\\[-8pt]\\[-8pt]
&&\qquad\lambda^d\bigl(\{f  \leq f(x)\} \cap B(x,s)\bigr)\geq c_f
\lambda^d(B(x,s)). \nonumber
\end{eqnarray}

\item[(b)] A Borel measurable map $f \dvtx \R^d \to\R$ has the \textit{convex sublevel approximation property} (CSAP) outside $\overline
{B}(0,R)$, $R>0$, if there is a~bounded convex set $C \subset\R^d$
with nonempty interior such that
\begin{eqnarray*}
&&\forall x \in\R^d \setminus\overline{B}(0,R), \exists
\varphi_x\dvtx \R^d \to\R^d,  \mbox{Euclidean motion}, \exists
a_x \geq1\\
&&\qquad \mbox{such that } x \in\varphi_x(a_x C) \subset\{f
\leq f(x)\}.
\end{eqnarray*}
[By Euclidean motion, we mean an affine transform of the
form $\varphi(y) =Ay+b$, $A$~orthogonal matrix and $b \in\R^d$.]\vspace*{2pt}

\item[(c)] A probability distribution $P$ has the \textit{peakless
sublevel tail property} (PSTP) outside $\overline B(0,R)$, $R>0$,
if:\vspace*{2pt}

\begin{enumerate}[(iii)]
\item[(i)]
$P$ is absolutely continuous with an essentially
bounded density $h$,\vspace*{2pt}
\item[(ii)]
$h$ is bounded away from $0$ on compacts sets, that is,
%
%
\begin{equation}\label{eq50}
\forall\rho>0, \exists c_{\rho} > 0 \qquad\mbox{such that } h(x)
\geq c_{\rho} \qquad\mbox{for all }
x \in\overline{B}(0,\rho).
\end{equation}

\item[(iii)]
There exist a~function $f\dvtx\R^d\to I$, $I$ interval of $\R$,
having the PSP and a~nonincreasing function $g\dvtx I\to(0,+\infty)$ such that
%
\[
\forall x \in\R^d\qquad \|x\| \ge R\quad\Longrightarrow\quad h(x)=g\circ f(x).
\]
Note that $\operatorname{supp}(P)=\R^d$.
\end{enumerate}
\end{enumerate}
\end{Dfn}
\begin{Pro}\label{42} If $f \dvtx \R^d \to\R^d$ has the CSAP outside
$\overline{B}(0,R)$, then it has the PSP outside $\overline{B}(0,R)$.
\end{Pro}
\begin{pf}
Let $s_0 > 0$ be arbitrary. By~\cite{Graf-Lu},
Example 12.7, there exists a~constant $\widetilde{c} > 0$ such that
%
%
\begin{equation}\label{eq48}
\forall x \in C, \forall s \in(0,s_0) \qquad \lambda
^d\bigl(C \cap B_{\| \ccccdot\|_2}(x,s)\bigr) \geq\widetilde{c} \lambda^d\bigl(B_{\|\ccccdot\|_2}(x,s)\bigr).
\end{equation}
%
There exists a~constant $\kappa\in(0,\infty)$ such that
\[
\frac{1}{\kappa} \| \cdot\|_2 \leq\| \cdot\| \leq\kappa
\| \cdot\|_2.
\]
Now let $x \in\R^d$ with $\|x\| \geq R$ and let $s \in
(0,s_0)$ be arbitrary. Then we have
\begin{eqnarray*}
&&\lambda^d\bigl(\{f \leq f(x)\} \cap B(x,s)\bigr) \\
&&\qquad \geq \lambda
^d \biggl(\varphi_x(a_x C) \cap B_{\|\ccccdot\|_2} \biggl(x, \frac
{s}{\kappa} \biggr) \biggr)\\
&&\qquad = \lambda^d \biggl(a_x C \cap\varphi_x^{-1}
\biggl(B_{\|\ccccdot\|_2} \biggl(x, \frac{s}{\kappa} \biggr) \biggr) \biggr)\\
&&\qquad = a_x^d \lambda^d \biggl(C \cap\frac{1}{a_x}
\varphi_x^{-1} \biggl(B_{\|\ccccdot\|_2}\biggl (x, \frac{s}{\kappa}
\biggr) \biggr) \biggr)\\
&&\qquad = a_x^d \lambda^d \biggl(C \cap B_{\|\ccccdot\|_2}
\biggl(\frac{1}{a_x} \varphi_x^{-1}(x), \frac{s}{a_x \kappa} \biggr)
\biggr)\\
&&\qquad \geq \widetilde{c} a_x^d \lambda^d \biggl(B_{\|\ccccdot\|
_2} \biggl(\frac{1}{a_x} \varphi_x^{-1}(x), \frac{s}{a_x
\kappa} \biggr) \biggr)\qquad \mbox{owing to (\ref{eq48})}\\
&&\qquad = \widetilde{c} a_x^d \frac{1}{\kappa^d a_x^d}
s^d \lambda^d \bigl(B_{\|\ccccdot\|_2}(0,1) \bigr)\\
&&\qquad = \widetilde{c} \kappa^{-d}
\frac{\lambda^d(B_{\|\ccccdot\|_2}(0,1))}{\lambda^d(B(0,1))} \lambda^d(B(x,s)).
\end{eqnarray*}
\upqed\end{pf}
\begin{Examples*}
\begin{enumerate}[(a)]
\item[(a)]
If \mbox{$\| \cdot\|_0$} is any norm on $\R^d$ and $f \dvtx\R^d
\to\R$ is defined by $f(x) = \|x\|_0$. Then $f$ has the CSAP outside
$\overline{B}(0,R)$, for every $R > 0$.

In particular, every nonsingular normal distribution has the PSTP
outside $\overline{B}(0,R)$ for every $R > 0$ and more generally, this
is the case for hyper-exponential distributions of the forms
\[
h(x)= K\|x\|_2^ae^{-c\|x\|_2^b},\qquad a,b,c, K>0.
\]
for large enough $R>0$ (in fact this is true for any norm).
\end{enumerate}
\end{Examples*}
\begin{pf}
Let $R > 0$ be arbitrary. Then there is an
$\widetilde{R} > 0$ with
\[
\overline{B}_{\|\ccccdot\|_0}(0, \widetilde{R}) \subset\overline{B}(0,R).
\]
Let $C = \overline{B}_{\|\ccccdot\|_0}(0, \widetilde{R})$.
Then $C$ is convex with nonempty interior. Let $x \in\R^d
\setminus\overline{B}_{\|\ccccdot\|_0}(0,\widetilde{R})$ be arbitrary.
Set $\varphi_x = id_{\R^d}$ and $a_x = \frac{1}{\widetilde{R}} \|
x\|_0 \geq1$. Then
\[
x = \varphi_x \biggl(a_x \widetilde{R} \frac{x}{\|x\|_0}
\biggr) \in\varphi_x(a_xC) = \overline{B}_{\|\ccccdot\|_0}
(0, \|x\|_0 ) = \{f \leq f(x)\}.
\]
\upqed\end{pf}

\begin{enumerate}[(b)]
\item[(b)]
\textit{Let $f \dvtx\R^d \to\R$ be semi-concave outside $\overline
{B}(0,R)$ in the following sense}:
\[
\exists\theta>1, \exists L > 0, \exists\varrho\dvtx\R^d
\setminus\overline{B}(0,R) \to\R_+ \setminus\{0\}, \exists
\delta\dvtx\R^d \setminus\overline{B}(0,R) \to\R^d \setminus\{0\}
\]
\textit{such that}:
\begin{enumerate}[(ii)]
\item[(i)] $\forall x \in\R^d \setminus\overline{B}(0,R)
, \frac{\varrho(x)}{\|\delta(x)\|_2} \leq L$,
\item[(ii)] $\forall x \in\R^d \setminus\overline{B}(0,R),
\forall y \in B (x, (\frac{1}{L} )^{
{1}/({\theta- 1})} )$,
$f(y) \leq f(x) + \delta(x) \cdot(y-x) + \varrho(x) \|y-x\|
_2^\theta$,
\textit{where $w \cdot z$ denotes the standard scalar product of} $w$, $z \in\R^d$.
\end{enumerate}
\textit{Then $f$ has the CSAP outside $\overline{B}(0,R)$.}
\end{enumerate}
\begin{pf}
Set $C = \{y=(y_1, \ldots, y_d) \in\R^d
\vert y_1+ L \|y\|_2^\theta\leq0\}$.
We will show that $C$ is a~bounded convex set with nonempty interior.
For $\lambda\in[0,1]$ and $y, \widetilde{y} \in C$ we have
\begin{eqnarray*}
&&\bigl(\lambda y_1 + (1- \lambda) \widetilde{y}_1 \bigr) + L \|
\lambda y + (1- \lambda) \widetilde{y}\|_2^\theta
\\
&&\qquad\leq\lambda y_1+(1- \lambda) \widetilde{y}_1 +
L\bigl(\lambda\|y\|_2 + (1- \lambda) \|\widetilde{y}\|_2\bigr)^\theta.
\end{eqnarray*}
Since $\theta> 1$, we have
\[
\bigl(\lambda\|y\|_2 + (1- \lambda) \|\widetilde{y}\|_2
\bigr)^\theta\leq\lambda\|y\|_2^\theta+ (1- \lambda) \|\widetilde
{y}\|_2^\theta,
\]
which yields
\[
\lambda y + (1- \lambda)\widetilde{y} \in C.
\]
Thus, $C$ is convex. For $y \in C$, we have
\begin{eqnarray*}
0 &\geq& y_1+L \|y\|_2^\theta \geq- \|y\|_2+ L \|y\|_2^\theta\\
&=& \|y\|_2 (L \|y\|_2^{\theta-1} - 1),
\end{eqnarray*}
hence $\|y\|_2 \leq(\frac{1}{L} )^{{1}/({\theta
-1})}$, so that $C$ is bounded.

There exists a~$t > 0$ with $-t + L t^\theta= t(L t^{\theta-1}-1) <
0$. For $y = (-t, 0, \ldots, 0)$ this implies $y_1+L \|y\|_2^\theta
< 0$. Hence, there exists a~neighborhood of $y$ which is contained in
$C$, that is, the interior of $C$ is not empty.

Now let $x \in\R^d$ with $\|x\| > R$ be arbitrary. Set $u = \frac
{\delta(x)}{\|\delta(x)\|_2}$. Let\vspace*{1pt} $\psi_x$ be a~rotation which maps
$e_1=(1,0, \ldots, 0)$ onto $u$. Define $\varphi_x \dvtx \R^d \to\R^d$
by $\varphi_x(y) = \psi_x(y) + x$. Then $\varphi_x$ is a~Euclidean
motion. Set $a_x = 1$. Since $0 \in C$ we have $x \in\varphi_x(C)
= \varphi_x(a_x C)$. For $y \in\varphi_x(a_x C) = \varphi_x(C)$
there is a~$z \in C$ with $y = \varphi_x(z)$, hence
%
\begin{eqnarray*}
\delta(x) \cdot(y-x) + \varrho(x) \|y-x\|_2^\theta &=& \delta(x)
\cdot\psi_x(z) + \varrho(x) \|\psi_x(z)\|_2^\theta\\
&=& \|\delta(x)\|_2 u \cdot\psi_x(z) + \varrho(x) \|\psi_x(z)\|
_2^\theta\\
&=& \|\delta(x)\|_2 e_1 \cdot z + \varrho(x) \|z\|_2^\theta\\
&=& \|\delta(x)\|_2 \biggl(z_1 + \frac{\varrho(x)}{\|\delta(x)\|
_2} \|z\|^\theta_2 \biggr)\\
&\leq& \|\delta(x)\|_2 (z_1+ L \|z\|_2^\theta) \leq0
\end{eqnarray*}
since $z \in C$. Moreover, $\|\varphi_x(z)-x\|_2 = \|\psi
_x(z)\|_2 = \|z\|_2$ and
$- \|z\|_2+ L \|z\|_2^\theta\leq0$ implies $\|z\|_2 \leq
(\frac{1}{L} )^{{1}/({\theta-1})}$, that is, $y= \psi_x(z)
\in B (x, ( \frac{1}{L} )^{{1}/({\theta-1})})$.

By (ii), this yields
\[
f(y) \leq f(x) + \delta(x) \cdot(y-x) + \varrho(x) \|y-x\|
_2^\theta\leq f(x)
\]
and, hence,
\[
\varphi_x(a_x C) \subseteq\{f \leq f(x)\}.
\]
\upqed\end{pf}

\begin{enumerate}[(c)]
\item[(c)]
\textit{Let $f \dvtx\R^d \to\R$ be a~differentiable function and let $R>0$
be such that there exist real constants $\alpha\in(0,1)$, $\beta>
0$ and $c \in(0,+ \infty)$ satisfying}:\vspace*{2pt}

\begin{enumerate}[(ii)]
\item[(i)]
$\forall x, y \in\R^d , [x,y] := \{x+t(y-x) , t \in
[0,1]\} \subset\R^d \setminus\overline{B}(0, R)
\Longrightarrow\break\|{\operatorname{grad} f(x)} - \operatorname{grad } f(y)\|
\leq c \|x-y\|^\alpha(1+ \|x\|^\beta+ \|y\|^\beta).$\vspace*{2pt}
\item[(ii)]
$ \inf_{\|x\| \geq R} \frac{\|
\operatorname{grad } f(x)\|}{1+ \|x\|^\beta} > 0$.
\end{enumerate}
\textit{Then $f$ is semi-concave outside of $\overline{B}(0,R + 1)$.}
\end{enumerate}
\begin{pf}
For every $x, y \in\R^d$ with $\|x\|>R$ and
$\|x-y\| \leq1$, we have
\[
\|y\|^\beta\leq(\|x\| + \|y-x\|)^\beta\leq(\|x\|
+1)^\beta= \|x\|^\beta\biggl(1+ \frac{1}{\|x\|} \biggr)^\beta
\]
%
so that
\begin{eqnarray*}
1+ \|x\|^\beta+ \|y\|^\beta &\leq& 1+ \|x\|^\beta\biggl(\biggl(1+ \frac
{1}{R} \biggr)^\beta+ 1 \biggr)\\
&\leq& \biggl( \biggl(1+ \frac{1}{R} \biggr)^\beta+ 1 \biggr) (\|x\|
^\beta+ 1 ).
\end{eqnarray*}
Let $\kappa\in(0,\infty)$ such that $\frac{1}{\kappa}\| \cdot\|_2
\leq\| \cdot\| \leq\kappa\| \cdot\|_2$.

Let $\theta=1+ \alpha$. Define $\varrho\dvtx \R^d \to\R_+ \setminus\{
0\}$ by $\varrho(x) = \kappa^2 c((1+ \frac{1}{R}
)^\beta+ 1 ) (\|x\|^\beta+1 )$ and $\delta\dvtx \R^d \to
\R^d$ by $\delta(x) = \operatorname{grad} f(x)$. Since $M :=
\inf_{\| x\| \geq R} \frac{\|{\operatorname{grad} f(x)}\|}{1+
\|
x\|^\beta} > 0$, we have $\delta(x) \neq0$ for all $x \in\R^d
\setminus\overline{B}(0,R)$. Moreover,
\[
\frac{\varrho(x)}{\|\delta(x)\|_2} \leq\frac{\varrho
(x)}{({1}/{\kappa}) \|\delta(x)\|} \leq\kappa^3 c \biggl(
\biggl(1+ \frac{1}{R} \biggr)^\beta+ 1 \biggr) \frac{1}{M} \leq L,
\]
where $L = \max\{1, \kappa^3 c ( (1+ \frac{1}{R}
)^\beta+ 1 ) \frac{1}{M} \}$. Let $x \in\R^d
\setminus\overline{B}(0,R+1)$ and $y \in B (x, (\frac
{1}{L} )^{{1}/({\theta-1})} )$ be arbitrary. Since $L
\geq1$ we have $[x,y] \subset\R^d \setminus\overline{B}(0, R)$
and, by the mean value theorem of differentiation,
\begin{eqnarray*}
f(y) - f(x) &=& ( \operatorname{grad } f(x) ) \cdot(y-x)\\
&&{} + \bigl( \operatorname{grad } f\bigl(x+t(y-x) \bigr) - \operatorname{grad } f(x) \bigr)
\cdot(y-x)
\end{eqnarray*}
for some $t \in[0,1]$. By our assumption, we obtain
\begin{eqnarray*}
&&\bigl(\operatorname{grad } f\bigl( x+t(y-x)\bigr) - \operatorname{grad } f(x)\bigr) \cdot(y-x)\\
&&\qquad \leq\bigl\|{\operatorname{grad} f}\bigl(x+t(y-x)\bigr) - \operatorname{grad } f(x)\bigr\|_2 \|y-x\|
_2\\
&&\qquad \leq\kappa^2 \bigl\|{\operatorname{grad} f}\bigl(x+t(y-x)\bigr) - \operatorname{grad } f(x)\bigr\|
\|y-x\|\\
&&\qquad \leq\kappa^2 ct^\alpha\|y-x\|^\alpha\bigl(1+ \|x\|^\beta+ \|
x+t(y-x)\|^\beta\bigr) \|y-x\|.
\end{eqnarray*}
%
Since $\|x+t(x-y) - x\| = t \|x-y\| \leq(\frac{1}{L}
)^{{1}/({\theta-1})} \leq1$, we deduce
\begin{eqnarray*}
&&\bigl(\operatorname{grad } f\bigl(x+t(y -x)\bigr) - \operatorname{grad } f(x)\bigr) \cdot(y-x)\\
&&\qquad \leq\kappa^2 c \biggl( \biggl(1+ \frac{1}{R} \biggr)^\beta+ 1
\biggr) (\|x\|^\beta+ 1 ) \|y-x\|^\theta\\
&&\qquad \leq\varrho(x) \|y-x\|^\theta.
\end{eqnarray*}
It follows that
\[
f(y) \leq f(x) + \delta(x) \cdot(y-x) + \varrho(x) \|y-x\|^\theta.
\]
Thus, $f$ is semi-concave outside the ball $\overline{B}(0,R+1)$.
\end{pf}

As always in this manuscript $\alpha_n$ is an $n$-optimal codebook for
$P$ of order $r > 0$, where we assume $\int\|x\|^{r+ \delta} \,dP(x)
< \infty$ for some $\delta> 0$.

Our first aim is to prove another variant of the first micro--macro
inequality for distributions $P$ having the PSTP.
\begin{Pro}\label{44}
Let $P$, with density $h$, have the PSTP outside $\overline B(0,R)$ for
a given $R>0$.
There exists a~constant $c_{21} > 0$ such that
%
%
\begin{eqnarray}\label{eq51}
&&\forall K \subset\R^d, \mbox{compact}, \exists n_K \in
\N \mbox{ such that } \forall n \geq n_K, \forall x
\in K \nonumber\\[-8pt]\\[-8pt]
&&\qquad c_{21} n^{-1/d} h(x)^{-{1}/({r+d})}  \geq
d(x, \alpha_n).\nonumber
\end{eqnarray}
%
%
\end{Pro}
\begin{pf}
Let $K \subset\R^d$ be compact. Since $\operatorname{supp}(P) = \R^d$, Proposition 2.2 in~\cite{Delattre} implies
\[
\lim_{n \to\infty} \max_{y \in K} d(y,
\alpha_n) = 0.
\]
Let $f$ and $g$ be as in Definition~\ref{41}(c)(iii) and let $s_0 >
0$ be related to $f$ by Definition~\ref{41}(a). Choose $n_K \in\N
$, so that
\[
\forall n \geq n_K \qquad \max_{y \in K} d(y, \alpha
_n) < \min(s_0,R).
\]
%
%
Let $n \geq n_K$ and let $x \in K$ be arbitrary. By (\ref{eq17}), we
know that
%
%
\begin{equation}\label{eq52}
c_5 (e^r_{n,r} - e^r_{n+1,r} ) \geq d(x, \alpha
_n)^{r+d} \frac{P(B(x,bd(x, \alpha_n)))}{\lambda^d(B(x, bd(x,
\alpha_n)))}.
\end{equation}
Since $\overline{B}(0,2R)$ is bounded and convex, there exists a~constant $\widetilde{c} > 0$ with
\[
\forall s \in(0,s_0), \forall y \in\overline{B}(0,2R)
\qquad \lambda^d\bigl(\overline{B}(0,2R) \cap B(y,s)\bigr) \geq\widetilde
{c} \lambda^d(B(y,s)).\vadjust{\goodbreak}
\]
If $x \in\overline{B}(0,2R)$, by Definition~\ref{41}(c)(ii)
there exists a~lower bound $c_{2R}>0$ of $h$ on $\overline{B}(0,2R)$,
so that
\begin{eqnarray*}
P(B(x,bd(x,\alpha_n))) & \ge & c_{2R} \lambda^d\bigl(\overline
{B}(0,R) \cap B(x,bd(x, \alpha_n))\bigr)\\
& \geq &c_{2R} \widetilde{c} \lambda^d(B(x,bd(x, \alpha_n))),
\end{eqnarray*}
hence $ c_5 (e^r_{n,r}-e^r_{n+1,r} ) \geq c_{2R}
\widetilde{c} d(x, \alpha_n)^{r+d}$
%
and consequently
%
%
\begin{equation}\label{eq53}
c_5 (e^r_{n,r} - e^r_{n+1,r} ) \geq c_{2R} \widetilde
{c} \frac{1}{\|h\|_{\overline{B}(0,2R)}} h(x) d(x, \alpha_n)^{r+d}
\end{equation}
for every $x \in\overline B(0$, $2R)$.
If $x \notin\overline{B}(0,2R)$ and $y \in B(x,bd(x, \alpha_n))
\cap\{f \leq f(x)\}$, then we have
\[
y\notin\overline B(0,R) \quad\mbox{and}\quad h(y) = g(f(y)) \geq
g(f(x)) = h(x)
\]
since $g$ is nonincreasing and we obtain
\begin{eqnarray*}
P(B(x,bd(x, \alpha_n)) )& \geq&P\bigl(B(x,bd(x, \alpha_n))
\cap\{f \leq f(x)\}\bigr)\\
& = & \int_{\{f \leq f(x)\} \cap B(x, bd(x, \alpha_n))} h(y) \,d \lambda^d(y)\\
& \geq&h(x) \lambda^d\bigl(\{f \leq f(x)\} \cap B(x, bd(x,\alpha
_n))\bigr)\\
& \geq& c_f h(x) \lambda^d(B(x,bd(x, \alpha_n)))
\end{eqnarray*}
since $f $ has the PSP. Hence,
%
%
\begin{equation}\label{eq54}
c_5 (e^r_{n,r} - e^r_{n+1,r} ) \geq c_f h(x) d(x, \alpha
_n)^{r+d}.
\end{equation}
Note that, by Proposition~\ref{25}, there exists a~constant $c_{11} >
0$ such that
\[
\forall n\in\N\qquad e^r_{n,r} -e^r_{n+1,r} \le c_{11} n^{- (1+r/d)}.
\]
Setting $c_{21} = ( c_{11} c_5 \max\{c^{-1}_f ,
(c_{2R} \widetilde{c} )^{-1} \} )^{
{1}/({r+d})}$ and combining the last inequality with (\ref{eq53}) and
(\ref{eq54}) yields the conclusion of the proposition.
\end{pf}
%
%
%
\begin{Remark*}
Note at this stage that the results established in the
rest of this section depend only on properties (\ref{eq50}) and (\ref
{eq51}), not directly on PSP.

Our next aim is to give an upper and a~lower bound for $P(W(a \vert
\alpha_n))$ and the local quantization error $\int_{W(a |
\alpha_n)}\|x-a\|^r\,dP(x)$, provided all the $W(a \vert\alpha
_n)$ intersect a~given compact set. The following lemma provides an
essential tool for the proof. Here and in the rest of the paper, we set
\[
\overline{s}_{n,a} = \sup\{\|x-a\| , x \in W(a \vert
\alpha_n)\},
\]
which can be considered as the \textit{radius} of the Voronoi cell $W(a
\vert\alpha_n)$.
\end{Remark*}
%
\begin{Lem}\label{46} Let $K \subset
\cop{\operatorname{supp}(P)}$ be an arbitrary compact set and let $\varepsilon
> 0$ be arbitrary. Then there exists an $n_{K, \varepsilon} \in\N$
such that
%
%
\begin{equation}\label{eq55}
\forall n \geq n_{K, \varepsilon}, \forall a \in\alpha_n
\qquad W(a \vert\alpha_n) \cap K \neq\varnothing
\quad\Rightarrow\quad\overline{s}_{n,a} \leq\varepsilon.\vadjust{\goodbreak}
\end{equation}
\end{Lem}

\begin{pf}
Let $\varepsilon> 0$. Since $K\subset\cop
{\operatorname{supp}(P)}$, one may assume without loss of generality that
$\varepsilon$ is small enough so that the $\varepsilon$-neighborhood
$K_{\varepsilon}:=\{y \in\R^d | d(y,K) \le\varepsilon\}$ is
included in $\operatorname{supp} P$. Since $K$ is compact and contained in
$\operatorname{supp}(P)$,~\cite{Delattre}, Proposition 2.2 implies $\lim
_{n \to\infty} \max_{x \in K} d(x, \alpha_n) = 0$.
Hence, there exists an $n_0 \in\N$ with
%
%
\begin{equation}\label{eq560}
\forall x \in K, \forall n \geq n_0 \qquad d(x, \alpha_n)
< \frac{\varepsilon}{2}.
\end{equation}
Now assume that (\ref{eq55}) does not hold for $\frac{\varepsilon
}{2}$ in the place of $\varepsilon$. Then there exist sequences
$(n_k)_{k \in\N}$ in $\N$ and $(a_k)$ with $n_k \uparrow\infty,
a_k \in\alpha_{n_k}$ with
\[
W(a_k \vert\alpha_{n_k}) \cap K \neq\varnothing,
\]
and $\overline{s}_{n_k,a_k} > \frac{\varepsilon}{2}$. Without loss
of generality, we assume $n_k > n_0$ for all $k \in\N$. For each $k
\in\N$, there is an $\widetilde{x}_k \in W(a_k, \alpha
_{n_k} )$ with $\|\widetilde{x}_k- a_k\| > \frac
{\varepsilon}{2}$. Set $x_k = a_k + \frac{\varepsilon}{2 \|
\widetilde{x}_k-a_k\|} (\widetilde{x}_k - a_k )$. Then
we have\vspace*{1pt} $\|x_k - a_k\| = \frac{\varepsilon}{2}$ and, since $W
(a_k, \alpha_{n_k} )$ is star shaped with center $a_k$
(see~\cite{Graf-Lu}, Proposition 1.2), we deduce that
$x_k \in[a_k, \widetilde{x}_k] \subset W(a_k \vert
\alpha_{n_k})$. Now let $z_k \in W(a_k \vert
\alpha_{n_k} ) \cap K$. Then $\|z_k -a_k\|< \frac{\varepsilon
}{2}$ owing to (\ref{eq560}) and $\|x_k-a_k\| = \frac{\varepsilon
}{2}$, so that $x_k \in K_\varepsilon$.

Since $K_\varepsilon$ is compact there exists a convergent
subsequence of $(x_k)$, whose limit we denote by $x_\infty\in
K_\varepsilon$. Then we have
\begin{eqnarray*}
d(x_\infty, \alpha_{n_k} ) & \geq & d(x_k, \alpha
_{n_k} ) - \|x_k -x_\infty\|\\
& = & \|x_k - a_k\| - \|x_k- x_\infty\|\\
& = & \frac{\varepsilon}{2} - \|x_k - x_\infty\|
\end{eqnarray*}
so that $\lim\sup_{k \to\infty}
d(x_\infty, \alpha_{n_k} ) \geq\frac{\varepsilon}{2}$.

Since $x_\infty\in K_\varepsilon\subset\operatorname{supp}(P)$, we know
that $\lim_{n \to\infty} d(x_\infty, \alpha_n) = 0$
(see~\cite{Graf-Lu}, Lemma 6.1 and~\cite{Delattre}, Proposition 2.2)
and obtain a~contradiction.
\end{pf}
\begin{Dfn} For a~compact set $K \subset\R^d$, let
\[
\alpha_n(K) = \{a \in\alpha_n \vert W(a \vert\alpha_n)
\cap K \neq\varnothing\}.
\]
\end{Dfn}
%
%
\begin{Pro}\label{48} Let $P$ satisfy the micro--macro inequality (\ref{eq51}).
There are constants $c_{22}, c_{23}, c_{24}, c_{25} > 0$ such that,
for every compact set $K \subset\R^d$ and every $\varepsilon>
0$, there exists an $n_{K, \varepsilon} \in\N$ such that, for
every $n \geq n_{K, \varepsilon}$, and every $a \in\alpha_n(K)$
the Voronoi cell $W(a \vert\alpha_n)$ is contained in
$K_\varepsilon$ and
%
%
\begin{eqnarray}\qquad
\label{eq56} P(W(a \vert\alpha_n)) & \leq&
c_{22} \bigl(\|h\|_{W(a \vert\alpha_n)} \bigr)^{
{r}/({r+d})} \frac{1}{n}, \\
\label{eq57} \int_{W(a \vert\alpha_n)} \|x-a\|
^r \,dP(x) & \leq& c_{23} \biggl( 1 + \log
\frac{\|h\|_{W(a \vert\alpha_n)}}{\operatorname{essinf}
h_{|W(a \vert\alpha_n)}} \biggr) n^{- (1+ {r/d})},
\\
\label{eq58} P (W_0(a \vert\alpha_n) ) &
\geq& c_{24} \bigl(\operatorname{essinf} h_{|W(a \vert
\alpha_n)} \bigr)^{{r}/({r+d})} \frac{1}{n} ,\\
\qquad
\label{eq59}
\int_{W_0(a \vert\alpha_n)} \|x-a\|^r \,dP(x)
& \geq&
c_{25} \biggl(\frac{\operatorname{essinf} h_{|W(a \vert\alpha
_n)} }{\|h\|_{W(a \vert\alpha_n)}} \biggr)^{\max(r,1)}
n^{- (1+ {r/d})}.
\end{eqnarray}
\end{Pro}
\begin{pf}
Let $K \subset\R^d$ be compact and
$\varepsilon> 0$ be arbitrary. By Lemma~\ref{46} and Proposition \ref
{44}, there exists an $n_{K, \varepsilon} \in\N$ with $n_{K,
\varepsilon} \geq2$ such that
%
%
\begin{equation}
\label{eq60} \forall n \geq n_{K, \varepsilon}, \forall a
\in\alpha_n(K) \qquad W(a \vert\alpha_n) \subset K_\varepsilon
\end{equation}
and
%
%
\begin{equation}\label{eq61} \forall n \geq n_{K, \varepsilon},
\forall x \in K_\varepsilon\qquad c_{21} n^{-1/d} h(x)^{-
{1}/({r+d})} \geq d(x, \alpha_n).
\end{equation}
Now let $n \geq n_{K, \varepsilon}$ and let $a \in\alpha_n(K)$ be
fixed. Set $\overline{t}_{n,a} = \|h\|_{W(a \vert\alpha_n)}$
and $\underline{t}_{n,a} = \operatorname{essinf} h_{|W(a \vert
\alpha_n)}$. Since $W(a \vert\alpha_n) \subset K_\varepsilon
$ by (\ref{eq60}), inequality (\ref{eq61}) implies
%
%
\begin{equation}\label{eq62}\qquad
\forall t > 0, \forall x \in\{h > t\} \cap W(a \vert
\alpha_n) \qquad \|x-a\| \leq c_{21}n^{-1/d} t^{- {1}/({r+d})}.
\end{equation}
This yields
%
%
\begin{eqnarray}
\label{eq63} \lambda^d\bigl(\{h > t\} \cap W(a \vert\alpha_n)\bigr)
&\leq& \lambda^d \bigl(B\bigl(a, c_{21} n^{-1/d} t^{-{1}/({r+d})}
\bigr)\bigr) \nonumber\\[-8pt]\\[-8pt]
& = & \lambda^d(B(0,1)) c^d_{21} t^{- {d}/({r+d})} n^{-1}.\nonumber
\end{eqnarray}
Now we will prove (\ref{eq56}).
Observing that $\lambda^d(\{h > t\} \cap W(a \vert\alpha
_n)) = 0$ for $t > \overline{t}_{n,a}$ we deduce
\begin{eqnarray*}
P(W(a \vert\alpha_n)) & =& \int_{W(a \vert
\alpha_n)} h \,d \lambda^d\\
& =& \int^\infty_0 \lambda^d\bigl(\{h>t\} \cap W(a \vert\alpha
_n)\bigr)\, dt\\
& =& \int^{\overline{t}_{n,a}}_0 \lambda^d\bigl(\{h > t\}
\cap W(a \vert\alpha_n)\bigr) \,dt\\
& \le& \biggl( \int_0^{\overline{t}_{n,a}} t^{- {d}/({r+d})}
\,dt \biggr) \lambda^d(B(0,1)) c^d_{21} n^{-1} \qquad\mbox{owing
to (\ref{eq63})}\\
& \leq& \lambda^d(B(0,1)) \frac{r+d}{r} c^d_{21} \bigl(\|h\|
_{W(a \vert\alpha_n)} \bigr)^{{r}/({r+d})} \frac{1}{n},
\end{eqnarray*}
which proves (\ref{eq56}) with $c_{22} = \lambda^d(B(0,1))\frac
{r+d}{r} c^d_{21}$.

Next, we will show (\ref{eq57}).
Using again $\lambda^d(\{h > t\} \cap W(a \vert\alpha
_n)) = 0$ for $t > \overline{t}_{n,a}$, we get
%
%
\begin{eqnarray}
\label{eq64}
\int_{W(a \vert\alpha_n)} \|x -a\|^r \,dP(x) & = & \int
_{W(a \vert\alpha_n)} \|x-a\|^r h(x) \,d \lambda^d(x)\nonumber\\
& = &\int_0^\infty\int_{\{h>t\} \cap W(a \vert\alpha_n)}
\|x-a\|^r \,d \lambda^d(x) \,dt\\
& = &\int_0^{\overline{t}_{n,a}} \int_{\{h>t\} \cap W(a \vert
\alpha_n)} \|x-a\|^r \,d \lambda^d(x) \,dt.\nonumber
\end{eqnarray}
For $t \leq\underline{t}_{n,a}$, we have $h(y) \geq t$ for $\lambda
^d\mbox{-a.e. } y \in W(a \vert\alpha_n)$ so that
\[
\int_{\{h>t\} \cap W(a \vert\alpha_n)} \|x-a\|
^r \,d \lambda^d(x) = \int_{W(a \vert\alpha_n)} \|x-a\|^r \,d \lambda^d(x).
\]
By (\ref{eq60}) and (\ref{eq61}), we have, for $\lambda^d\mbox
{-a.e. } x \in W(a \vert\alpha_n)$,
\[
\|x-a\| = d(x, \alpha_n) \leq c_{21} n^{-1/d} h(x)^{-{1}/({r+d})} \leq c_{21} n^{-1/d} (\underline{t}_{n,a}
)^{-{1}/({r+d})}
\]
so that
\[
\lambda^d\bigl(W(a \vert\alpha_n)\setminus B\bigl(a, c_{21}
n^{-1/d} (\underline{t}_{n,a} )^{-{1}/({r+d})}
\bigr) \bigr)=0.
\]
%
Consequently,
%
%
\begin{eqnarray}\label{eq66}\qquad
&&\int_0^{\underline{t}_{n,a}} \int_{\{h>t\} \cap
W(a \vert\alpha_n)} \|x-a\|^r \,d \lambda^d(x)\,
dt\nonumber\\
&&\qquad \leq \int^{\underline{t}_{n,a}}_0 \int_{_{B (a, c_{21}
n^{-1/d}(\underline t_{n,a})^{-{1}/({r+d})})}}
\bigl(c_{21} n^{-1/d}(\underline t_{n,a})^{-{1}/({r+d})}\bigr)^r \,d \lambda^d(x) \,dt\\
&&\qquad = c_{23} n^{-(1+r/d)},\nonumber
\end{eqnarray}
%
%
%
where $c_{23} = c_{21}^{r+d} \lambda^d(B(0,1))
$.
Using (\ref{eq62}) and the same argument as before, we obtain
%
%
\begin{eqnarray}
\label{eq67}
&&\int^{\overline{t}_{n,a}}_{\underline{t}_{n,a}} \int
_{\{h>t\} \cap W(a \vert\alpha_n)} \|x-a\|^r \,d \lambda^d(x) \,dt \nonumber\\
&&\qquad \leq\int_{\underline{t}_{n,a}}^{\overline{t}_{n,a}}
\int_{B(a, c_{21} n^{-1/d} t^{-{1}/({r+d})})}
c_{21}^rt^{-{r}/({r+d})}n^{-r/d} \,dP(x) \,dt\nonumber\\[-8pt]\\[-8pt]
&&\qquad \leq c_{23}n^{- (1+ {r/d})} \int_{\underline
{t}_{n,a}}^{\overline{t}_{n,a}} t^{-1} \,dt\nonumber\\
&&\qquad = c_{23} n^{- (1+ {r/d} )} \log
\biggl( \frac{\overline{t}_{n,a}}{\underline{t}_{n,a}}\biggr).\nonumber
\end{eqnarray}
Combining (\ref{eq66}) and (\ref{eq67}) with (\ref{eq64})
yields (\ref{eq57}).

Now we will prove (\ref{eq58}).
It follows from the second micro--macro inequality (Proposition \ref
{23}) and Proposition~\ref{25} that there exists a~real constant $c >
0$ (independent of $n$ and $a$) such that
%
%
\begin{equation}\label{eq68}
cn^{- (1+ {r/d} )} \leq\int_{W_0(a \vert
\alpha_n)} \bigl(d(x, \alpha_n \setminus\{a\})^r - \|x-a\|^r\bigr) \,dP(x).
\end{equation}
%
%
Since (4.59) implies that $\ovli{W_0(a \vert\alpha_n)}$ is
compact and nonempty there exists a~$z \in\partial W_0(a \vert
\alpha_n)$. Obviously this $z$ satisfies
\[
\|z-a\| = d(z,\alpha_n \setminus\{a\})
\]
and, therefore,
%
%
\begin{equation}\label{eq69}d(a, \alpha_n \setminus\{a\}) \leq\|
a-z\| + d(z, \alpha_n \setminus\{a\}) = 2 \|z-a\|.
\end{equation}
This implies that, for every $x \in\ovli{W_0(a \vert\alpha_n)}$,
\begin{eqnarray*}
d(x,\alpha_n \setminus\{a\}) & \leq&\|x-a\| + d(a, \alpha_n
\setminus\{a\})\\
& \leq&\|x-a\| + 2 \|z-a\| = d(x, \alpha_n) + 2d(z, \alpha_n).
\end{eqnarray*}
Since $d_{\alpha_n} := d(\cdot, \alpha_n)$ is continuous and every
nonempty relatively open subset of $\ovli{W_0(a \vert\alpha
_n)}$ has positive Lebesgue measure, we deduce
\[
\max\{d(y, \alpha_n)\dvtx y \in\ovli{W_0(a \vert\alpha_n)}\} =
\operatorname{esssup} d_{\alpha_n \vert\ovli{W_0(a \vert\alpha_n)}}.
\]
By (\ref{eq60}) and (\ref{eq61}) this yields
\begin{eqnarray*}
d(x, \alpha_n \setminus\{a\}) & \leq & 3 \operatorname{esssup} d_{\alpha
_n \vert\ovli{W_0(a \vert\alpha_n)}}\\
& \leq & 3c_{21} n^{-1/d} \operatorname{esssup} \bigl(h_{\vert\ovli
{W_0(a \vert\alpha_n)}} \bigr)^{-{1}/({r+d})}\\
& = & 3c_{21} n^{-1/d} \bigl(\operatorname{essinf} h_{\vert\ovli{W_0(a
\vert\alpha_n)}} \bigr)^{-{1}/({r+d})}\\
& \leq & 3c_{21} n^{-1/d} (\underline{t}_{n,a} )^{- {1}/({r+a})}
\end{eqnarray*}
and, therefore,
%
%
\begin{equation}\label{eq71}
\int_{W_0(a \vert\alpha_n)}\! d(x, \alpha_n
\setminus\{a\})^r \,dP(x) \leq3^r c_{21}^r n^{-r/d}
(\underline t_{n,a} )^{- {r}/({r+d})} P(W_0(a \vert
\alpha_n)).\hspace*{-35pt}
\end{equation}
Using (\ref{eq68}), we deduce
\[
c 3^{-r} c_{21}^{-r} (\underline{t}_{n,a} )^{
{r}/({r+d})} n^{-1} \leq P(W_0(a \vert\alpha_n))
\]
and, hence, (\ref{eq58}) with $c_{24} = c 3^{-r} c^{-r}_{21}$.

Now we will prove (\ref{eq59}). It follows from (\ref{eq68}) that
%
%
\begin{equation}\label{eq72}
c n^{- (1+ {r/d} )} \leq\int_{W_0(a \vert
\alpha_n)}\! \bigl(\bigl(\|x-a\|+ d(a, \alpha_n
\setminus\{a\}) \bigr)^r - \|x-a\|^r\bigr) \,dP(x).\hspace*{-35pt}
\end{equation}

\setcounter{Case}{0}
\begin{Case}[($r \geq1$)] Using the mean value theorem for
differentiation yields
%
%
\begin{eqnarray}\label{eq73}
&&c n^{- (1+ {r/d} )}\nonumber\hspace*{-35pt}\\[-8pt]\\[-8pt]
&&\qquad\leq\int_{W_0(a \vert\alpha_n)} r \bigl( \|
x-a\| + d(a, \alpha_n \setminus\{a\}) \bigr)^{r-1}
d(a, \alpha_n
\setminus\{a\}) \,dP(x).\nonumber\hspace*{-35pt}
\end{eqnarray}
By (\ref{eq69}), (\ref{eq60}) and (\ref{eq61}), we know that
%
%
\begin{equation}\label{eq74}
\|x-a\| + d(a, \alpha_n \setminus\{a\}) \leq3c_{21} n^{-1/d}
(\underline{t}_{n,a} )^{-{1}/({r+d})}.
\end{equation}
Combining (\ref{eq73}) and (\ref{eq74}) yields
%
%
\begin{equation}\label{eq75}
cn^{- (1+{r/d} )} \leq d(a, \alpha_n \setminus\{
a\}) r \bigl(3 c_{21} n^{-1/d} (\underline{t}_{n,a}
)^{-{1}/({r+d})}\bigr)^{r-1} P(W_0(a \vert\alpha_n)).\hspace*{-35pt}
\end{equation}
By (\ref{eq56}), we have
\[
P(W_0(a \vert\alpha_n)) \leq c_{22} \overline{t}_{n,a}^{
{r}/({r+d})} \frac{1}{n}
\]
and, hence,
%
%
\begin{equation}\label{eq76}\qquad
c_{22}^{-1} cr^{-1}(3c_{21})^{1-r} \underline{t}_{n,a}^{
({r-1})/({r+d})} \overline{t}_{n,a}^{-{r}/({r+d})} n^{-1/d} \leq
d(a, \alpha_n \setminus\{a\}).
\end{equation}
Set $\widetilde{c} = c_{22}^{-1} cr^{-1}(3c_{21})^{1-r}$.
Then we deduce
%
%
\begin{equation}\label{eq77}
B \biggl(a, \frac{\widetilde{c} }{2} \underline{t}_{n,a}^{
({r-1})/({r+d})} \overline{t}_{n,a}^{- {r}/({r+d})} n^{-1/d}
\biggr) \subset W_0(a \vert\alpha_n).
\end{equation}

It follows that
%
%
\begin{eqnarray}\label{eq78}
&&\int_{B (a, ({ \widetilde{c}}/{2 }) \underline
{t}_{n,a}^{({r-1})/({r+d})} \overline{t}_{n,a}^{- {r}/({r+d})}
n^{-1/d})} \|x-a\|^r h(x) \,d \lambda^d(x) \nonumber\\[-8pt]\\[-8pt]
&&\qquad\leq
\int_{W_0(a \vert\alpha_n)} \|x-a\|^r \,dP(x).\nonumber
\end{eqnarray}
Since $h(x) \geq\underline{t}_{n,a}$, for $\lambda^d\mbox
{-a.e. } x \in B(a, \frac{\widetilde{c}}{2} \underline
{t}_{n,a}^{({r-1})/({r+d})} \overline{t}_{n,a}^{- {r}/({r+d})}
n^{-1/d})$ and
\[
\int_{B(a,\varrho)} \|x-a\|^r \,d \lambda^d(x) = \varrho^{r+d} \int
_{B(0,1)} \|u\|^r \,d \lambda^d(u)
\]
for every $\varrho> 0$, the left-hand side of (\ref{eq78}) is greater
or equal to
\begin{eqnarray*}
&&\underline{t}_{n,a} \int_{B(0,1)} \|x\|^r \,d \lambda^d(x)
\biggl(\frac{\widetilde{c} }{2} \underline{t}_{n,a}^{({r-1})/({r+d})}
\overline{t}_{n,a}^{- {r}/({r+d})} \biggr)^{r+d} n^{- (1+ {r/d} )}\\
&&\qquad = \int_{B(0,1)} \|u\|^r \,d \lambda^d(u)\biggl(\frac{\widetilde
{c} }{2} \biggr)^{r+d} \underline{t}_{n,a}^r \overline
{t}_{n,a}^{-r} n^{- (1+ r/d )}.
\end{eqnarray*}
Inequality (\ref{eq59}) follows by setting $c_{25} = \int
_{B(0,1)} \|u\|^r \,d \lambda^d(u) (\frac{\widetilde{c}
}{2} )^{r+d}$.
\end{Case}
\begin{Case}[($r < 1$)] In this case, we have
\[
\bigl(\|x-a\| + d(a, \alpha_n \setminus\{a\})\bigr)^r \leq\|x-a\|^r+ d(a,
\alpha_n \setminus\{a\})^r\vadjust{\goodbreak}
\]
for all $x \in W_0(a \vert\alpha_n)$, so that, by (\ref{eq72}),
%
%
\begin{eqnarray}
\label{eq79} cn^{-(1+{r/d} )} & \leq &\int
_{W_0(a \vert\alpha_n)} d(a, \alpha_n \setminus\{a\})^r
\,dP(x) \nonumber\\[-9pt]\\[-9pt]
& \leq & d(a, \alpha_n \setminus\{a\})^r P(W_0(a \vert
\alpha_n)).\nonumber
\end{eqnarray}
By (\ref{eq56}), we know that
\[
P(W_0(a \vert\alpha_n)) \leq c_{22} (\overline{t}_{n,a}
)^{{r}/({r+d})} \frac{1}{n}
\]
and, hence,
%
%
\begin{equation}\label{eq80}
c^{1/r} c_{22}^{-{1}/{r}} \overline{t}_{n,a}^{ -
{1}/({r+d})} n^{-1/d} \leq d(a, \alpha_n \setminus\{a\}).
\end{equation}
As above this implies, for $\widetilde{c} = c^{1/r} c_{22}^{-1/r}$,
\[
\underline{t}_{n,a} \int_{B(0,1)} \|x\|^r \,d \lambda^d(x) \biggl(
\frac{ \widetilde{c}}{2 } \biggr)^{r+d} \frac{\underline
{t}_{n,a}}{\overline{t}_{n,a}} n^{- (1+ {r/d} )}
\leq\int_{W_0(a \vert\alpha_n)} \|x-a\|^r \,dP(x)
\]
and (\ref{eq59}) follows.\qed\vspace*{-3pt}
\end{Case}
\noqed\end{pf}
\begin{Thm}\label{49}
Let $P$ satisfy the micro--macro inequality (\ref{eq51}). Then there
are constants $c_{22}, c_{23}, c_{24}, c_{25} > 0$ such that, for every
compact set $K \subset\R^d$, the following hold:
%
%
\begin{eqnarray}
\label{eq81}
&\displaystyle \limsup_{n \to\infty} n \max_{a
\in\alpha_n(K)} P(W(a \vert\alpha_n)) \leq c_{22}
\Bigl( \inf_{\varepsilon> 0} \|h\|_{K_\varepsilon}
\Bigr)^{{r}/({r+d})},&
\\[-2pt]
%
%
\label{eq82}
&\displaystyle \limsup_{n \to\infty} n^{1+r/d}
\max_{a \in\alpha_n(K)}
\int_{W(a \vert\alpha_n)} \|x-a\|^r \,dP(x)&\nonumber\\[-9pt]\\[-9pt]
&\displaystyle \qquad\leq c_{23} \biggl(1 + \log\biggl(\inf_{\varepsilon>
0} \frac{\|h\|_{K_\varepsilon}}{\operatorname{essinf}
h_{|K_{\varepsilon}} }\biggr) \biggr),\qquad\qquad&\nonumber
%
%
\\[-2pt]
\label{eq83}
&\displaystyle \liminf_{n \to\infty} n
\min_{a \in\alpha_n(K)} P(W_0(a \vert
\alpha_n)) \geq c_{24} \sup_{\varepsilon> 0}
(\operatorname{essinf} h_{|K_\varepsilon} )^{{r}/({r+d})},&
\\[-2pt]
%
%
\label{eq84}
&\displaystyle \liminf_{n \to\infty} n^{(1+
{r/d} )} \min_{a \in\alpha
_n(K)} \int_{W(a \vert\alpha_n)}
\|x-a\|^r \,dP(x) &\nonumber\\[-9pt]\\[-9pt]
&\displaystyle \geq c_{25} \sup_{\varepsilon>
0} \biggl( \frac{\operatorname{essinf} h_{|K_\varepsilon}}{\|h\|
_{K_{\varepsilon}} } \biggr)^{\max(1,r)}.\qquad\quad&\nonumber\vspace*{-3pt}
\end{eqnarray}
\end{Thm}
\begin{pf}
The theorem follows immediately from
Proposition~\ref{48}.\vspace*{-3pt}
\end{pf}
\begin{Cor} \label{410} For every $x \in\R^d$, let $a_{n,x} \in
\alpha_n$ satisfy $x \in W(a_{n,x} \vert\alpha_n)$. Then
%
%
\begin{eqnarray}
\label{eq85} &\displaystyle \limsup_{n \to\infty} nP(W(a_{n,x}
\vert\alpha_n)) \leq c_{22} \Bigl( \limsup_{y \to x}
h(y) \Bigr)^{{r}/({r+d})},&
\\[-2pt]
%
%
\label{eq86}
&\displaystyle \limsup_{n \to\infty} n^{1+r/d}\int_{W(a_{n,x} \vert\alpha_n)}
\|x-a\|^r \,dP(x)& \nonumber\\[-8pt]\\[-8pt]
&\displaystyle \quad\hspace*{5pt}\leq c_{23} \biggl(1+ \log\lim
_{\varepsilon\downarrow0} \frac{\sup h(B(x,\varepsilon))}{
\inf h(B(x, \varepsilon))}\biggr),&\nonumber
\\
%
%
\label{eq87}
&\displaystyle \liminf_{n \to\infty} nP(W_0(a_{n,x},
\vert\alpha_n)) \geq c_{24} \Bigl(\liminf_{y \to
x} h(y) \Bigr)^{{r}/({r+d})},&
\\
%
%
\label{eq88}
&\displaystyle \liminf_{n \to\infty} n^{1+r/d}
\int_{W_0(a_{n,x} \vert\alpha_n}
\|x-a\|^r \,dP(X) &\nonumber\\[-8pt]\\[-8pt]
&\displaystyle \quad\geq c_{25} \biggl(\lim_{\varepsilon
\downarrow0} \frac{ \inf h(B(x, \varepsilon))}{\sup h(B(x,
\varepsilon))} \biggr)^{\max(1,r)}.&\nonumber
\end{eqnarray}
Moreover, if $h$ is continuous, then
$\limsup_{y \to x} h(y) = h(x) = \liminf
_{y \to x} h(y)$ and
\[
\lim_{\varepsilon\downarrow0} \frac{\sup h(B(x,
\varepsilon))}{ \inf h(B(x, \varepsilon)} = \lim
_{\varepsilon\downarrow0} \frac{ \inf h(B(x, \varepsilon
))}{\sup h(B(x, \varepsilon))}= 1 .
\]
\end{Cor}
\begin{pf}
The corollary follows from Theorem 4.9 if one sets $K= \{x\}$.~%
\end{pf}
\begin{Remarks*}
(a) For certain one-dimensional distribution
functions, sharper versions of the above corollary have been proved by
Fort and Pag\`{e}s (\cite{Fort}, Theorem 6).

(b) If $R > 0$ and the density $h$ has the form $h(x) = g(\|
x\|_0)$ for all $x \notin B(0,R)$, where $g\dvtx [0,+\infty) \to
(0,+\infty)$ is a~decreasing function and $\| \cdot\|_0$ is an
arbitrary norm on $\R^d$ then there exists a~constant $c > 0$ and an
$m = m(c) \in\N$ such that
\[
\forall n \geq m, \forall x \in\R^d \qquad c n^{-1/d}
h(x)^{-{1}/({r+d})} \geq d(x,\alpha_n).
\]
This can be used to show that there is a~$\widetilde{c} > 0$ with
\[
\forall n \geq m, \forall a \in\alpha_n\qquad P(W(a \vert\alpha_n)) \leq
\widetilde{c} \bigl(\|h\|_{W(a \vert\alpha_n)} \bigr)^{
{r}/({r+d})} \frac{1}{n}.
\]
Under additional assumptions on $g$ ($g$ regularly varying), one can
also give a~similar upper bound for the local $L^s$-quantization
errors, $s \in(0,r)$.
\end{Remarks*}

\section{The local quantization behavior in the interior of the
support}\label{sec5}

In this section, we will show that weaker versions of the results in
Section~\ref{4} still hold without assuming the strong version of the first
micro--macro inequality as stated in (\ref{eq51}). We have to restrict
our investigations to compact sets in the interior of the support of
the probability in question and also obtain weaker constants in the
corresponding inequalities for the local probabilities and quantization errors.

Let $r \in(0, \infty)$ be fixed. In this section $P$ is always an
absolutely continuous Borel probability\vadjust{\goodbreak} on $\R^d$ with density $h$. We
assume that there is a~$\delta> 0$ with $ \int\|x\|^{r+ \delta
} \,dP(x) < +\infty$. As before, $\alpha_n$ is an $n$-optimal
codebook for $P$ of order $r$. For $n \in\N$ and $a \in\alpha
_n$ set $\overline{s}_{n,a} = \sup\{\|x-a\| , x \in W(a
\vert\alpha_n)\}$ and $\underline{s}_{n,a} = \sup\{s > 0 ,
B(a, s) \subset W(a \vert\alpha_n)\}$.

Moreover,\vspace*{-3pt} we assume that $h$ is essentially bounded and
that $ \operatorname{essinf} h_{|K} > 0$ for every compact set $K
\subset\cop{\operatorname{supp}(P)}$, where\vspace*{2pt} $\stackover{B}$ denotes the
interior of the set $B \subset\R^d$. For the use in the first
micro--macro inequality, we fix a~$b \in(0, \frac{1}{2} )$.
%
\begin{Lem} \label{51} There exists a~constant $c_{26} > 0$ such that,
for every $n \in\N$ and $a \in\alpha_n$,
%
%
\begin{equation}
\label{eq89} c_{26} n^{-1/d} \bigl(\operatorname{essinf} h_{|B(a,
(1+b) \overline{s}_{n,a})} \bigr)^{-{1}/({r+d})} \geq
\overline{s}_{n,a}.
\end{equation}
\end{Lem}
\begin{pf}
By the first micro--macro inequality (\ref
{eq17}) and Proposition~\ref{25} there exists a~constant $c > 0$ with
%
%
\begin{equation} \label{eq90}\qquad
\forall n \in\N, \forall x \in\R^d\qquad cn^{-
(1 + r/d)} \geq d(x,\alpha_n)^{r+d} \frac{P(B(x,bd(x, \alpha
_n)))}{\lambda^d(B(x, bd(x, \alpha_n))}.
\end{equation}

Now let $n \in\N$ and $a \in\alpha_n$ be arbitrary.

It follows from (\ref{eq90}) that
%
%
\begin{equation} \label{eq91}\qquad
\forall x \in W(a \vert\alpha_n) \qquad \|x-a\|^{r+d}
\frac{P(B(x,b\|x-a\|))}{\lambda^d(B(x, b\|x-a\|))} \leq cn^{-
(1 + {r/d} )}.
\end{equation}
For $x \in W(a \vert\alpha_n)$ and $y \in B(x,
bd(x, \alpha_n))$, we have
\[
\|y-a\| < \|y-x\| + \|x-a\| \leq b \|x-a\| + \|x-a\| \leq(1+b)
\overline{s}_{n,a}
\]
so that
%
%
\begin{equation}
\label{eq92} B(x,b\|x-a\|) \subseteq B\bigl(a, (1+b) \overline{s}_{n,a}\bigr).
\end{equation}
This yields
%
%
\begin{eqnarray}
\label{eq93} P\bigl(B(x,b\|x-a\|)\bigr) & = & \int_{B(x,b\|x-a\|)}
h\,d \lambda^d \nonumber\\[-8pt]\\[-8pt]
&\ge& \operatorname{essinf} h_{|B(a, (1+b) \overline
{s}_{n,a})} \lambda^d\bigl(B(x,b\|x-a\|)\bigr)\nonumber
\end{eqnarray}
owing to (\ref{eq92}). Thus, (\ref{eq91}) implies
%
%
\begin{equation}
\label{eq4} \|x-a\|^{r+d} \operatorname{essinf} h_{|B(a,(1+b)
\overline{s}_{n,a})} \leq c n^{- (1+ {r/d})}.
\end{equation}
Since $x \in W(a \vert\alpha_n)$ was arbitrary,
we deduce
\[
\overline{s}_{n,a}^{r+d} \operatorname{essinf} h_{|B(a,(1+b)
\overline{s}_{n,a})} \leq c n^{- (1+ {r/d} )}
\]
and, hence, (\ref{eq89}) with $c_{26} = c^{{1}/({r+d})}$.
\end{pf}
%
\begin{Lem} \label{52}There exist real constants $c_{27}, c_{28}>0$
such that, for every $n \in\N$ and $a \in\alpha_n$,
%
%
\begin{equation}
\label{eq95}
P(W(a \vert\alpha_n)) \leq c_{27} \frac{\|h\|
_{B(a, \overline{s}_{n,a})}}{(\operatorname{essinf} h_{|B(a,(1+b)
\overline{s}_{n,a})})^{{d}/({r+d})}} n^{-1}\vadjust{\goodbreak}
\end{equation}
and
\begin{equation}
\label{eq96}
\int_{W(a \vert\alpha_n)} \|x-a\|^r \,dP(x)
\leq c_{28} \frac{\|h\|_{B (a, \overline{s}_{n,a})}}{
\operatorname{essinf} h_{|B(a(1+b), \overline{s}_{n,a})}} n^{-(1+
{r/d} )}.
\end{equation}
\end{Lem}
\begin{pf}
Let $n \in\N$ and $a \in\alpha_n$ be arbitrary. Then (\ref
{eq89}) implies
\begin{eqnarray*}
P(W(a \vert\alpha_n)) & \leq & P(B(a, \overline{s}_{n,a})) \leq\|
h\|_{B(a, \overline{s}_{n,a})} \lambda^d(B(a, \overline
{s}_{n,a}))\\
& \leq &\lambda^d(B(0,1)) \|h\|_{B(a, \overline{s}_{n,a})} \overline
{s}_{n,a}^d\\
& \le &\lambda^d(B(0,1)) c_{26}^d \|h\|_{B(a, \overline{s}_{n,a})}
\bigl( \operatorname{essinf} h_{|B(a,(1+b) \overline{s}_{n,a})}
\bigr)^{- {d}/({r+d})} n^{-1}.
\end{eqnarray*}
Thus (\ref{eq95}) follows for $c_{27} = \lambda^d(B(0,1))
c_{26}^d$.

Similarly (\ref{eq89}) implies
\begin{eqnarray*}
&&\int_{W(a \vert\alpha_n)} \|x-a\|^r \,dP(x)
\\
&&\qquad \leq \int
_{B(a, \overline{s}_{n,a})} \|x-a\|^r \,dP(x)\\
&&\qquad \leq \|h\|_{B(a, \overline{s}_{n,a})} \int_{B(a, \overline
{s}_{n,a})} \|x-a\|^r \,d \lambda^d(x)\\
&&\qquad \leq \lambda^d(B(0,1)) \|h\|_{B(a, \overline{s}_{n,a})}
\overline{s}_{n,a}^{r+d}\\
&&\qquad \le \lambda^d(B(0,1))c_{26}^{r+d} \|h\|_{B(a, \overline
{s}_{n,a})} \bigl(\operatorname{essinf} h_{|B(a,(1+b) \overline
{s}_{n,a})} \bigr)^{-1} n^{- (1+ {r/d})}
\end{eqnarray*}
still owing to (\ref{eq89}).

Thus, (\ref{eq96}) follows for $c_{28} = \lambda^d(B(0,1))
c_{26}^{r+d}$.
\end{pf}
%
\begin{Lem} \label{53}There exists real constants $c_{29}, c_{30}> 0$
such that, for every $n \geq2$ and every $a \in\alpha_n$,
%
%
\begin{equation}
\label{eq97a} \underline{s}_{n,a}\geq c_{29} \frac{(
\operatorname{essinf} h_{|B(a,(1+b) \overline{s}_{n,a}} )^{1-
{1}/({r+d})}}{\|h\|_{B(a,\overline{s}_{n,a})}} n^{-1/d} \qquad\mbox{for
} r \geq1
\end{equation}
and
%
%
\begin{equation}
\label{eq97b}
\underline{s}_{n,a} \geq c_{30} \biggl( \frac{(
\operatorname{essinf} h_{|B(a,(1+b)
\overline{s}_{n,a})})^{{d}/({r+d})}}{\|h\|_{B(a, \overline{s}_{n,a})}}
\biggr)^{1/r} n^{-1/d} \qquad\mbox{for } 0 < r < 1.\hspace*{-40pt}
\end{equation}
\end{Lem}
\begin{pf}
By the second micro--macro inequality
(Proposition~\ref{23}) combined with Proposition~\ref{25}, there is a~constant $c > 0$ such that
\[
\forall n \geq2 \qquad c n^{- (1+r/d)} \leq\int_{W_0(a
\vert\alpha_n)} \bigl(d(x, \alpha_n \setminus\{a\})^r
- \|x-a\|^r\bigr) \,dP(x).
\]

\setcounter{Case}{0}
\begin{Case}[($r \geq1$)] As in (\ref{eq72}) and (\ref
{eq73}), we deduce
%
%
\begin{eqnarray}
\label{eq98}
&&cn^{- (1+ {r/d} )} \leq\int
_{W_0(a \vert\alpha_n)} r\bigl(\|x-a\| + d(a, \alpha
_n \setminus\{a\})\bigr)^{r-1}\nonumber\\[-9pt]\\[-9pt]
&&\qquad\hspace*{75.3pt}{}\times d(a, \alpha_n \setminus\{a\})
\,dP(x).\nonumber\vspace*{-2pt}
\end{eqnarray}
Since $n \geq2$ there exists an $\widetilde{a} \in\alpha
_n \setminus\{a\}$ with
\[
W(a \vert\alpha_n) \cap W(\widetilde{a} \vert\alpha_n)
\neq\varnothing.\vspace*{-2pt}
\]
Let $z \in W(a \vert\alpha_n) \cap W(\widetilde
{a} \vert\alpha_n)$ be arbitrary. Then we have
\[
\|z-a\| = d(z,\alpha_n) = \|z- \widetilde{a}\|\vspace*{-2pt}
\]
and, hence
\[
d(a, \alpha_n \setminus\{a\}) \leq\|a- \widetilde{a}\| \leq
\|a-z\| + \|z- \widetilde{a}\| = 2\|z-a\|\vspace*{-2pt}
\]
so that
\[
d(a, \alpha_n \setminus\{a\}) \leq2 \overline{s}_{n,a}.\vspace*{-2pt}
\]
It follows from (\ref{eq98}) that
%
%
\begin{eqnarray}
\label{eq99}
cn^{-(1 + {r/d} )} &
\leq& r(3 \overline{s}_{n,a})^{r-1} d(a, \alpha_n \setminus\{a\})
P(W_0(a \vert\alpha_n)) \nonumber\\[-2pt]
& \leq & r(3 \overline{s}_{n,a})^{r-1} d(a, \alpha_n
\setminus\{a\}) \|h\|_{B(a, \overline{s}_{n,a})} \lambda
^d(B(0,1)) \overline{s}^d_{n,a}\\[-2pt]
& = &r3^{r-1} \overline{s}_{n,a}^{r+d-1} \lambda
^d(B(0,1)) \|h\|_{B(a, \overline{s}_{n,a})} d(a, \alpha_n
\setminus\{a\}).\nonumber\vspace*{-2pt}
\end{eqnarray}
%
%
This implies
\[
c r^{-1}3^{1-r} (\lambda^d(B(0,1)) )^{-1} \bigl(\|h\|
_{B(a, \overline{s}_{n,a})} \bigr)^{-1} \overline
{s}_{n,a}^{1-(r+d)} n^{- (1+ {r/d} )} \leq d(a,
\alpha_n \setminus\{a\})\vspace*{-2pt}
\]
and, hence, by (\ref{eq89})
\begin{eqnarray*}
&&cr^{-1} 3^{1-r} (\lambda^d(B(0,1)) )^{-1} \bigl(\|h\|
_{B(a, \overline{s}_{n,a})} \bigr)^{-1} c_{26}^{1-(r+d)}\\[-2pt]
&&\quad{}\times
\bigl(\operatorname{essinf} h_{|B(a, (1+b) \overline{s}_{n,a})}
\bigr)^{- ({1-(r+d)})/({r+d})}n^{-1/d}
\\[-2pt]
&&{}\qquad\leq d(a, \alpha_n \setminus\{a\}).\vspace*{-2pt}
\end{eqnarray*}
Since $\underline{s}_{n,a} = \frac{1}{2} d(a, \alpha_n
\setminus\{a\})$ this leads to (\ref{eq97a}) with
\[
c_{29} = \tfrac{1}{2} cr^{-1} 3^{1-r} (\lambda
^d(B(0,1)))^{-1} c_{26}^{1-(r+d)}.\vspace*{-2pt}
\]
\end{Case}
\begin{Case}[($r \leq1$)] As in (\ref{eq79}), we have
\begin{eqnarray*}
cn^{-1+r/d} & \leq & d(a, \alpha_n \setminus\{a\})^r P(W_0(a \vert
\alpha_n))\\[-2pt]
& \leq & d(a, \alpha_n \setminus\{a\})^r \|h\|_{B(a, \overline
{s}_{n,a})} \lambda^d(B(0,1)) \overline{s}_{n,a}^d\vspace*{-2pt}
\end{eqnarray*}
and, hence, by (\ref{eq89})
\begin{eqnarray*}
&&c n^{-( 1+ {r/d} )} \bigl(\|h\|_{B(a,
\overline{s}_{n,a})} \bigr)^{-1} (\lambda^d(B(0,1))
)^{-1} c_{26}^{-d} n \bigl(\operatorname{essinf} h_{|B(a,(1+b)
\overline{s}_{n,a})} \bigr)^{{d}/({r+d})} \\[-2pt]
&&\qquad\leq d(a, \alpha_n
\setminus\{a\})^r,\vspace*{-2pt}\vadjust{\goodbreak}
\end{eqnarray*}
which implies
\begin{eqnarray*}
&&c^{1/r} \bigl(\|h\|_{B(a, \overline{s}_{n,a})}
\bigr)^{-{1/r}} (\lambda^d(B(0,1) ))^{-1/r}
c_{26}^{- {d/r}} \\
&&\quad{}\times\bigl(\operatorname{essinf} h_{|B(a,(1+b)
\overline{s}_{n,a})} \bigr)^{{d}/({r(r+d)})} n^{-1/d} \\
&&{}\qquad\leq d(a,
\alpha_n \setminus\{a\}).
\end{eqnarray*}
Since $\underline{s}_{n,a} = \frac{1}{2} d(a, \alpha_n
\setminus\{a\})$ this leads to
\[
c_{30} \biggl( \frac{(\operatorname{essinf} h_{b(a,(1+b)
\overline{s}_{n,a})} )^{{d}/({r+d})}}{\|h\|_{B(a, \overline
{s}_{n,a})}} \biggr)^{1/r} n^{-1/d} \leq\underline{s}_{n,a}
\]
with
\[
c_{30} = \tfrac{1}{2} c^{1/r} (\lambda
^d(B(0,1)))^{-1/r} c_{26}^{- {d/r}}.
\]
\end{Case}
\upqed\end{pf}
%
%
\begin{Lem}\label{54} There exist constants $c_{31}, c_{32}, c_{33},
c_{34} > 0$ such that, for every $n > 2$ and $a \in\alpha_n$,
%
%
\begin{eqnarray}
\label{eq100}\qquad
&&P(W_0(a \vert\alpha_n)) \nonumber\\[-4pt]\\[-8pt]
&&\qquad\geq
\cases{
\displaystyle c_{31} \biggl( \frac{\operatorname{essinf} h_{|B(a,(1+b) \overline
{s}_{n,a})}}{\|h\|_{B(a,\overline{s}_{n,a})}} \biggr)^d
\bigl(\operatorname{essinf} h_{B(a,(1+b) \overline{s}_{n,a})}
\bigr)^{{r}/({r+d})} n^{-1},\vspace*{2pt}\cr
\qquad\mbox{for } r \geq1,\vspace*{2pt}\cr
\displaystyle c_{32} \biggl( \frac{\operatorname{essinf} h_{|B(a,(1+b) \overline
{s}_{n,a})}}{\|h\|_{B(a, \overline{s}_{n,a})}} \biggr)^{{d/r}}
\bigl(\operatorname{essinf} h_{|B(a,(1+b) \underline{s}_{n,a})}
\bigr)^{{r}/({r+d})} n^{-1},\vspace*{2pt}\cr
\qquad\mbox{for } 0 < r < 1,}\hspace*{-10pt}\nonumber
\end{eqnarray}
and
%
%
\begin{eqnarray}
\label{eq101}\qquad
&&\int_{W_0(a \vert\alpha_n)}
\|x-a\|^r \,dP(x)\nonumber\\[-4pt]\\[-8pt]
&&\qquad\geq
\cases{
\displaystyle c_{33} \biggl(\frac{\operatorname{essinf} h_{|B(a,(1+b) \overline
{s}_{n,a})}}{\|h\|_{B(a, \overline{s}_{n,a})}} \biggr)^{r+d}
n^{-(1+ {r/d} )}, &\quad for $r \geq1$,\vspace*{2pt}\cr
\displaystyle c_{34} \biggl( \frac{\operatorname{essinf} h_{|B(a,(1+b) \overline
{s}_{n,a})}}{\|h\|_{B(a, \overline{s}_{n,a})}} \biggr)^{1+
{d/r}} n^{- (1+ {r/d} )}, &\quad for $0 < r
< 1$.}\hspace*{-10pt}\nonumber
\end{eqnarray}
\end{Lem}
\begin{pf}
First, we will prove (\ref{eq100}). We have
\begin{eqnarray*}
P(W_0(a \vert\alpha_n)) & \geq&P(B(a, \underline{s}_{n,a})) =
\int_{B(a,\underline{s}_{n,a})} h\,d \lambda^d\\
& \geq&\operatorname{essinf} h_{|B(a,\underline{s}_{n,a})} \lambda
^d(B(0,1)) \underline{s}^d_{n,a}\\
& \geq&\operatorname{essinf} h_{|B(a,(1+b) \overline{s}_{n,a})}
\lambda^{d}(B(0,1)) \underline{s}^d_{n,a}.
\end{eqnarray*}
Using (\ref{eq97a}), we obtain
\begin{eqnarray*}
P(W_0(a \vert\alpha_n)) &\geq&\lambda^d(B(0,1)) c^d_{29}
\biggl( \frac{\operatorname{essinf} h_{|B(a,(1+b) \overline
{s}_{n,a})}}{\|h\|_{B(a, \overline{s}_{n,a})}} \biggr)^d \\
&&{}\times\bigl(
\operatorname{essinf} h_{|B(a,(1+b) \underline{s}_{n,a})} \bigr)^{
{r}/({r+d})} n^{-1}
\end{eqnarray*}
for $r \geq1$ and using (\ref{eq97b}) we get
\begin{eqnarray*}
P(W_0(a \vert\alpha_n)) & \geq &\lambda^d(B(0,1)) c_{30}^d
\bigl(\|h\|_{B(a, \overline{s}_{n,a})} \bigr)^{-{d/r}}\\
&&{}\times
\bigl(\operatorname{essinf} h_{|B(a,(1+b) \overline{s}_{n,a})}
\bigr)^{({d}/({r+d}))({d/r})+1} n^{-1}\\
& = &\lambda^d(B(0,1)) c^d_{30} \biggl( \frac{\operatorname{essinf}
h_{|B(a,(1+b) \overline{s}_{n,a})}}{\|h\|_{B(a, \overline
{s}_{n,a})}} \biggr)^{{d/r}} \\
&&{}\times \bigl(\operatorname{essinf}
h_{|B(a,(1+b) \overline{s}_{n,a})}\bigr)^{{r}/({r+d})} n^{-1}
\end{eqnarray*}
for $0 < r < 1$.
With $c_{31} = \lambda^d(B(0,1)) c^d_{29}$ and $c_{32} = \lambda
^d(B(0,1)) c^d_{30}$ we deduce~(\ref{eq100}).

Now we will prove (\ref{eq101}). We have
\begin{eqnarray*}
\int_{W_0(a \vert\alpha_n)} \|x-a\|^r \,dP(x)
& \geq& \int_{B(a, \underline{s}_{n,a})} \|x-a\|^r
\operatorname{essinf} h_{|B(a, \underline{s}_{a,n})} \,d \lambda^d(x)\\
& \geq& \bigl(\operatorname{essinf} h_{|B(a, \underline{s}_{n,a})}
\bigr) \int_{B(a, \underline{s}_{n,a})} \|x-a\|^r
\,d \lambda^d(x).
\end{eqnarray*}
%
Now
\[
\int_{B(a, \underline{s}_{n,a})} \|x-a\|^r \,d \lambda^d(x) = \underline{s}_{n,a}^{r+d} \int_{B(0,1)} \|x\|^r
\,d\lambda^d(x)
\]
so that
\[
\int_{W_0(a \vert\alpha_n)} \|x-a\|^r \,dP(x)
\geq\int_{B(0,1)} \|x\|^r \,d\lambda^d(x)
\operatorname{essinf} h_{|B(a, \underline{s}_{n,a})} \underline
{s}^{r+d}_{n,a}.
\]
Using Lemma~\ref{53}, we obtain (\ref{eq101}) with $c_{33}
= \int_{B(0,1)} \|x\|^r \,d\lambda^d(x) c_{29}^{r+d}$ and $c_{34} =
\int_{B(0,1)} \|x\|^r \,d\lambda^d(x)
c_{30}^{r+d}$.
\end{pf}
\begin{Lem} \label{55} Let $K \subset\cop{\operatorname{supp}(P)}$ be an
arbitrary compact set and let
\[
\varepsilon\in\bigl(0,d\bigl(K, \R^d \setminus\cop{\operatorname{supp}(P)}\bigr)\bigr)
\]
be arbitrary [where $d(K, \varnothing) = \infty$]. Then there exists an
$n_{K, \varepsilon} \in\N$ such that
%
%
\begin{equation}
\label{eq102} \forall n \geq n_{K, \varepsilon}, \forall a
\in\alpha_n(K) \qquad \overline{s}_{n,a} \leq\varepsilon,
\end{equation}
where $\alpha_n(K) = \{a \in\alpha_n \vert W(a \vert
\alpha_n) \cap K \neq\varnothing\}$.
\end{Lem}
\begin{pf}
The proof is identical to that of Lemma \ref{46}.
\end{pf}
\begin{Thm} \label{56}Let $P$ be an absolutely continuous Borel
probability measure on $\R^d$ with density\vspace*{1pt} $h$ and $ \int\|x\|^{r+
\delta} \,dP(x) < \infty$ for some $\delta> 0$. Then there exist
constants $c_{27}, c_{28}, c_{31}, c_{32}, c_{33}, c_{34} > 0$ such
that, for every compact\vspace*{-2pt} $K \subset\cop{\operatorname{supp}(P)}$, the following hold:
%
%
\begin{eqnarray}
\label{eq103}
&&\limsup_{n \to\infty} n \max_{a \in\alpha_n(K)} P(W(a
\vert\alpha_n)) \nonumber\\[-8pt]\\[-8pt]
&&\qquad\leq c_{27}
\inf_{\varepsilon> 0} \frac{\|h\|_{K_\varepsilon
}}{( \operatorname{essinf} h_{K_\varepsilon} )^{{d}/({r+d})}},\nonumber
\\
%
%
\label{eq104}
&&\limsup_{n \to\infty} n^{1+r/d} \max_{a \in\alpha_n(K)} \int_{W(a \vert
\alpha_n)} \|x-a\|^r \,dP(x) \nonumber\\[-8pt]\\[-8pt]
&&\qquad\leq c_{28}
\inf_{\varepsilon> 0} \frac{\|h\|_{K_\varepsilon}}{
\operatorname{essinf}
h_{|K_\varepsilon}},\nonumber
\\
%
%
\label{eq105}
&&\liminf_{n \to\infty} n \min_{a \in\alpha
_n(K)} P(W_0(a \vert\alpha_n)) \nonumber\\[-8pt]\\[-8pt]
&&\qquad\geq
\cases{
\displaystyle c_{31} \inf_{\varepsilon> 0} \biggl(\frac{
\operatorname{essinf} h_{|K_\varepsilon}}{\|h\|_{K_\varepsilon}}
\biggr)^d (\operatorname{essinf} h_{|K_\varepsilon})^{{r}/({r+d})}, \vspace*{2pt}\cr
\qquad\mbox{for } r \geq1,\vspace*{2pt}\cr
\displaystyle c_{32} \inf_{\varepsilon> 0} \biggl(\frac{
\operatorname{essinf} h_{|K_\varepsilon}}{\|h\|_{K_\varepsilon}}
\biggr)^{{d/r}} (\operatorname{essinf} h_{|K_\varepsilon}
)^{{r}/({r+d})}, \vspace*{2pt}\cr
\qquad\mbox{for } 0 < r < 1,}\nonumber
\end{eqnarray}
and
%
%
\begin{eqnarray}
\label{eq106}
&&\liminf_{n \to\infty} n^{1+
{r/d}} \min_{a \in\alpha_n(K)} \int_{W_0(a
\vert\alpha_n)} \|x-a\|^r \,dP(x) \nonumber\\[-8pt]\\[-8pt]
&&\qquad\geq
\cases{\displaystyle
c_{33} \inf_{\varepsilon> 0} \biggl(\frac{
\operatorname{essinf} h_{|K_\varepsilon}}{\|h\|_{K_\varepsilon}}
\biggr)^{r+d}, \vspace*{2pt}\cr
\qquad\mbox{for } r \geq1,\vspace*{2pt}\cr
\displaystyle c_{34} \inf_{\varepsilon> 0} \biggl(\frac{
\operatorname{essinf} h_{|K_\varepsilon}}{\|h\|_{K_\varepsilon}}
\biggr)^{1+{d/r}}, \vspace*{2pt}\cr
\qquad\mbox{for } 0 < r < 1.}\nonumber
\end{eqnarray}
\end{Thm}
\begin{pf}
Let $\varepsilon> 0$ satisfy $\varepsilon<
d(K, \R^d \setminus\cop{\operatorname{supp}(P)})$.\vspace*{-2pt} By Lemma~\ref{55} there
exists an $n_{K, \varepsilon} \in\N$ such that
\[
\forall n \geq n_{K, \varepsilon} , \forall a \in\alpha_n(K)
\qquad \overline{s}_{n,a} < \frac{\varepsilon}{2(1+b)}.
\]
This implies
\[
\forall n \geq n_{K, \varepsilon} , \forall a \in\alpha_n(K)
\qquad B\bigl(a, (1+b) \overline{s}_{n,a}\bigr) \subset K_\varepsilon
\]
%
and, therefore,
\[
\|h\|_{B(a,(1+b) \overline{s}_{n,a})} \leq\|h\|_{K_\varepsilon}
\]
as well as
\[
\operatorname{essinf} h_{|B(a,(1+b) \overline{s}_{n,a})} \geq
\operatorname{essinf} h_{|K_\varepsilon}
\]
for all $n \geq n_{K, \varepsilon}$ and all $a \in\alpha_n(K)$.

These inequalities combined with Lemma~\ref{52} and Lemma~\ref{54}
yield the assertions of the theorem.
\end{pf}
\begin{Remark*}
The above theorem yields estimates for the
asymptotics of the local cell probabilities and quantization errors
only if the density $h$ is essentially bounded and bounded away from
$0$ on each compact subset of the interior of the support of $P$.
\end{Remark*}
\begin{Cor} \label{58} For every $x \in\R^d$, let $a_{n,x} \in
\alpha_n$ satisfy $x \in W(a_{n,x} \vert\alpha_n)$. Assume
that $x \in\cop{\operatorname{supp}(P)}$ and $h$ is continuous at $x$. Then
\begin{eqnarray}\qquad
\min(c_{31}, c_{32}) h(x)^{
{r}/({r+d})}
&\leq&\liminf_{n \to\infty}
nP(W_0(a_{n,x} \vert\alpha_n))\nonumber\\[-8pt]\\[-8pt]
&\leq& \limsup_{n \to\infty}
nP(W(a_{n,x} \vert\alpha_n))
\leq c_{27} h(x)^{
{r}/({r+d})}\nonumber
\end{eqnarray}
and
\begin{eqnarray}\quad
\min(c_{33}, c_{34}) &\leq& \liminf_{n \to\infty
} n^{1+r/d} \int_{W(a_{n,x} \vert
\alpha_n)} \|y - a_{n,x}\|^r \,dP(y)\nonumber\\[-8pt]\\[-8pt]
&\leq& \limsup_{n \to\infty}
n^{1+r/d} \int_{W(a_{n,x} \vert\alpha_n)} \|y- a_{n,x}\|^r
\,dP(y) \leq c_{28}.\nonumber
\end{eqnarray}
\end{Cor}
\begin{pf}
Set $K = \{x\}$ in Theorem~\ref{56}.
\end{pf}


%

%
\printaddresses

\end{document}